\documentclass[A4j,11pt]{article}
\usepackage{graphics}
\usepackage{amsmath}
\usepackage{amssymb}
\usepackage{latexsym}
\usepackage{amsfonts}
\begin{document}

\title{{\bf On curvature-adapted and\\
proper complex equifocal submanifolds}}
\author{{\bf Naoyuki Koike}}
%
%
\date{}
%
\maketitle
\begin{abstract}
In this paper, we investigate curvature-adapted and proper complex 
equifocal submanifolds in a symmetric space of non-compact type.  
The class of these submanifolds contains principal orbits of Hermann type 
actions as homogeneous examples and is included by that of curvature-adapted 
and isoparametric submanifolds with flat section.  
First we introduce the notion of a focal point of non-Euclidean type 
on the ideal boundary for 
a submanifold in a Hadamard manifold and describe the equivalent condition 
for a curvature-adapted and complex equifocal submanifold to be proper 
complex equifocal in terms of this notion.  
Next we show that the complex Coxeter group associated with 
a curvature-adapted and proper complex equifocal submanifold is the same type 
group as one associated with a principal orbit of a Hermann type action 
and evaluate from above the number of distinct principal curvatures of 
the submanifold.  
\end{abstract}

\vspace{0.5truecm}

\section{Introduction}
In symmetric spaces, the notion of an equifocal submanifold was introduced 
in [TT].  This notion is defined as a compact submanifold with globally flat 
and abelian normal bundle such that the focal radius functions for each 
parallel normal vector field are constant.  
However, the equifocality is rather weak property in the case 
where the symmetric spaces are of non-compact type and the submanifold is 
non-compact.  So we [Koi1,2] 
have recently introduced the notion of a complex equifocal 
submanifold in a symmetric space $G/K$ of non-compact type.  
This notion is defined by imposing the constancy of the complex focal 
radius functions instead of focal radius functions.  
Here we note that the complex focal radii are the 
quantities indicating the positions of the focal points of the 
extrinsic complexification of the submanifold, where the submanifold needs 
to be assumed to be complete and of class $C^{\omega}$ 
(i.e., real analytic).  
On the other hand, Heintze-Liu-Olmos [HLO] has recently defined the 
notion of an isoparametric submanifold with flat section in a general 
Riemannian manifold as a submanifold such that the normal holonomy group is 
trivial, its sufficiently close parallel submanifolds are of constant mean 
curvature with respect to the radial direction and that the image of the 
normal space at each point by the normal 
exponential map is flat and totally geodesic.  
We [Koi2] showed the following fact:

\vspace{0.3truecm}

{\sl All isoparametric submanifolds with 
flat section in a symmetric space $G/K$ of non-compact type are complex 
equifocal and that conversely, all curvature-adapted and complex equifocal 
submanifolds are isoparametric ones with flat section.}

\vspace{0.3truecm}

Here the curvature-adaptedness means that, for each normal vector $v$ of 
the submanifold, the Jacobi operator $R(\cdot,v)v$ preserves the tangent space 
of the submanifold invariantly and the restriction of $R(\cdot,v)v$ to the 
tangent space commutes with the shape operator $A_v$, where $R$ is the 
curvature tensor of $G/K$.  Note that curvature-adapted hypersurfaces in 
a complex hyperbolic space (and a complex projective space) mean so-called 
Hopf hypersurfaces and that curvature-adapted complex equifocal hypersurfaces 
in the space mean Hopf hypersurfaces with constant principal curvatures, 
whcih are classified by J. Berndt [B1].  Also, he [B2] classified 
curvature-adapted hypersurfaces with constant prinicipal curvatures 
(i.e., curvature-adapted complex equifocal hypersurfaces) 
in the quaternionic hyperbolic space.  
In Appendix 2, we will prove an important fact 
for a curvature-adapted submanifold.  
As a subclass of the class of complex equifocal submanifolds, 
we [Koi3] defined that of the proper complex equifocal submanifolds in $G/K$ 
as a complex equifocal submanifold whose lifted submanifold to 
$H^0([0,1],\mathfrak g)$ ($\mathfrak g:={\rm Lie}\,G$) through 
some pseudo-Riemannian submersion of $H^0([0,1],\mathfrak g)$ onto $G/K$ is 
proper complex isoparametric in the sense of [Koi1], where we note that 
$H^0([0,1],\mathfrak g)$ is a pseudo-Hilbert space.  
For a proper complex equifocal $C^{\omega}$-submanifold, 
we [Koi4] defined the notion of the associated complex Coxeter group 
as the Coxeter group generated by the complex reflections of 
order two with respect to complex focal hyperplanes in the normal space of 
the lift of the complexification of the submanifold to 
$H^0([0,1],\mathfrak g^{\bf c})$ ($\mathfrak g^{\bf c}:={\rm Lie}\,G^{\bf c}$) 
by some anti-Kaehler submersion of $H^0([0,1],\mathfrak g^{\bf c})$ onto 
the anti-Kaehler symmetric space $G^{\bf c}/K^{\bf c}$, where we note that 
the lifted submanifold is proper anti-Kaehler isoparametric in the sense of 
[Koi2].  Here we note that the associated complex 
Coxeter group can be defined in terms of complex focal radii of 
the original submanifold without the use of the lifted submanifold.  
We [Koi4] showed that a proper complex equifocal submanifold 
is decomposed into the (non-trivial) extrinsic product of such submanifolds 
if and only if the associated complex Coxeter group is decomposable.  
Thus it is worth to investigate the complex Coxeter group in detail.  
According to Theorem 1 of [Ch], all complete equifocal 
submanifolds of codimension 
greater than one on simply connected irreducible symmetric space of 
compact type are homogeneous.  Hence they are principal orbits of 
hyperpolar actions (see [HPTT]).  
According to the classification of the hyperpolar actions by A. Kollross 
([Kol]), all hyperpolar actions of cohomogeneity greater than one on the 
irreducible symmetric space are Hermann ones.  
On the other hand, O. Goertsches and G. Thorbergsson ([GT]) has recently 
showed that principal orbits of Hermann actions are curvature-adapted.  
Hence we have the following fact:

\vspace{0.3truecm}

{\sl All complete equifocal submanifolds of 
codimension greater than one in simply}

{\sl connected irreducible symmetric 
spaces of compact type are catched as principal}

{\sl orbits of Hermann actions and hence they are curvature-adapted.}

\vspace{0.3truecm}

Let $G/K$ be a symmetric space of non-compact type and $H$ be a symmetric 
subgroup of $G$ such that $({\rm Fix}\,\sigma)_0\subset H\subset
{\rm Fix}\,\sigma$ for some involution $\sigma$ of $G$, where 
${\rm Fix}\,\sigma$ is the fixed point group of $\sigma$ and 
$({\rm Fix}\,\sigma)_0$ is the identity component of the group.  
We ([Koi2]) called the action of such a group $H$ on $G/K$ an 
{\it action of Hermann type}.  In this paper, we call this action 
Hermann type action for simplicity.  
We ([Koi2,3]) showed the following fact:

\vspace{0.3truecm}

{\sl Principal orbits of a Hermann type action are curvature-adapted 
and proper 

complex equifocal.}

\vspace{0.3truecm}

\noindent
From these facts, it is conjectured that comparatively many ones among 
complex equifocal submanifolds of codimension 
greater than one in irreducible symmetric spaces of non-compact type are 
curvature-adapted and proper complex equifocal.  
The following questions are naturally proposed:

\vspace{0.2truecm}

\noindent
{\bf Question.} {\sl Do all curvature-adapted and proper complex equifocal 
submanifolds occur as principal orbits of Hermann type actions?}

\vspace{0.2truecm}

\noindent
We defined the notion of a proper complex equifocal submanifold as a complex 
equifocal submanifold whose lifted submanifold to the above path space is 
a proper complex isoparametric submanifold.  
It is important to give an equivalent condition for a complex equifocal 
submanifold to be proper complex equifocal by using geometric quantities 
of the original submanifold without the use of those of the lifted 
submanifold.  
In this paper, we give such an equivalent condition for a curvature-adapted 
and complex equifocal submanifold.  For its purpose, we first 
introduce the notion of a focal point of non-Euclidean type 
on the ideal boundary $N(\infty)$ 
for a submanifold in a Hadamard manifold $N$ in general.  
By using this notion, we obtain the following equivalent condition.  

\vspace{0.5truecm}

\noindent
{\bf Theorem A.} {\sl Let $M$ be a curvature-adapted and complex equifocal 
submanifold in a symmetric space $N:=G/K$ of non-compact 
type.  Then the following conditions ${\rm (i)}$ and ${\rm (ii)}$ are 
equivalent:

{\rm(i)} $M$ is proper complex equifocal, 

{\rm (ii)} $M$ has no focal point of non-Euclidean type on the ideal boundary 
$N(\infty)$.}

\vspace{0.5truecm}

According to this theorem, we can catch a curvature-adapted and proper complex 
equifocal submanifold as a curvature-adapted and isoparametric submanifold 
with flat section which has no focal point of non-Euclidean type on the ideal 
boundary.  
In Section 6 of [Koi4], we investigated the complex Coxeter groups 
associated with principal orbits of Hermann type actions.  According to 
the investigation and Appendix of this paper, it follows that the complex 
Coxeter group associated 
with a principal orbit $H(gK)$ of a Hermann type action $H\times G/K\to G/K$ 
is isomorphic to the affine Weyl group (which is denoted by 
$W^A_{\overline{\triangle}}$) 
associated with the root system 
$\overline{\triangle}:=
\{\alpha\vert_{g_{\ast}^{-1}T^{\perp}_{gK}(H(gK))}\,\vert\,\alpha\in\triangle
\,\,{\rm s.t.}\,\,\alpha\vert_{g_{\ast}^{-1}T^{\perp}_{gK}(H(gK))}\not=0\}$, 
where $\triangle$ is the root system of $G/K$ with respect to a maximal 
abelian subspace $\mathfrak a$ containing 
$g_{\ast}^{-1}T^{\perp}_{gK}(H(gK))$.  
See Section 2 about the definition of the affine Weyl group associated with 
a root system.  
In order to make sure of whether the above question is solved affirmatively, 
it is important to investigate whether the complex Coxeter group associated 
with a curvature-adapted and proper complex equifocal submanifold is 
isomorphic to the same type group.  For the complex Coxeter group associated 
with this submanifold, we have the following fact.  

\vspace{0.5truecm}

\noindent
{\bf Theorem B.} {\sl Let $M$ be a curvature-adapted and proper complex 
equifocal $C^{\omega}$-submanifold in a symmetric space $G/K$ of 
non-compact type and $\triangle$ be the root system of $G/K$ 
with respect to a maximal abelian subspace $\mathfrak a$ of $T_{eK}(G/K)$ 
containing $\mathfrak b:=g_{\ast}^{-1}T^{\perp}_{gK}M$, 
where $gK$ is an arbitrary point of $M$.  
Then $\overline{\triangle}:=\{\alpha\vert_{\mathfrak b}\,\,\vert\,\,\alpha
\in\triangle\,\,{\rm s.t.}\,\,\alpha\vert_{\mathfrak b}\not=0\}$ is 
a weakly root system and the complex Coxeter group associated with $M$ is 
isomorphic to the affine Weyl group associated with $\overline{\triangle}$.}

\vspace{0.5truecm}

See Section 2 about the definition of a weakly root system.  
Thus the complex Coxeter group associated with a curvature-adapted and proper 
complex equifocal $C^{\omega}$-submanifold is isomorphic to the same type one 
as the group associated with a principal orbit of a Hermann type action.  
Hence the possibility for Question 2 to be solved affirmatively goes up.  

\vspace{0.5truecm}

\noindent
{\it Remark 1.1.} According to this theorem, in case of 
${\rm codim}\,M=1$, the complex Coxeter group associated with $M$ 
is isomorphic to ${\bf Z}_2\propto{\bf Z}$.  

\vspace{0.5truecm}

By using Theorem 2 of [Koi4] and Theorem B, we obtain the following 
splitting theorem.  

\vspace{0.5truecm}

\noindent
{\bf Corollary B.1.} {\sl Let $M$ and 
$\overline{\triangle}$ be as in Theorem B.  
Then $M$ is decomposed into the extrinsic 
product of two curvature-adapted and proper complex equifocal submanifolds 
if and only if $W_{\overline{\triangle}}$ 
is decomposable, where $W_{\overline{\triangle}}$ is 
the Coxeter group associated with 
$\overline{\triangle}$.}

\vspace{0.5truecm}

See Section 2 about the definition of the Coxeter group associated with 
a weakly root system.  
From this corollary, the following fact follows directly.  

\vspace{0.5truecm}

\noindent
{\bf Corollary B.2.} {\sl Let $M$ be as in Theorem B.  
If $G/K$ is reducible and ${\rm codim}\,M={\rm rank}\,G/K$, 
then $M$ is decomposed into the extrinsic product of two curvature-adapted 
and proper complex equifocal submanifolds.}

\vspace{0.5truecm}

For the number of mutually distinct principal curvatures of 
a curvature-adapted and proper complex equifocal $C^{\omega}$-submanifold, 
we have the following fact.  

\vspace{0.5truecm}

\noindent
{\bf Theorem C.} {\sl Let $M$ be a curvature-adapted and proper complex 
equifocal $C^{\omega}$-submanifold in a symmetric space $G/K$ of non-compact 
type and $A$ be the shape tensor of $M$.  Then, for each normal vector $v$ of 
$M$ at $gK$, we have 
$$
\sharp\,{\rm Spec}\,A_v\leq
\sharp(\overline{\triangle}_+\setminus\overline{\triangle}^1_+)\times 2
+\sharp\overline{\triangle}^1_+
+{\rm dim}\,\mathfrak z_{\mathfrak p}(\mathfrak b)-{\rm codim}\,M,$$
where ${\rm Spec}\,A_v$ is the spectrum of $A_v$, $\overline{\triangle}$ is 
as in the statement of Theorem B, 
$\overline{\triangle}^1_+:=\{\beta\in\overline{\triangle}_+\,\vert\,
{\rm the}\,\,
{\rm multiplicity}\,\,{\rm of}\,\,\beta\,\,{\rm is}\,\,{\rm equal}\,\,
{\rm to}\,\,1\}$, 
$\sharp(\cdot)$ is the cardinal number of $(\cdot)$ and 
$\mathfrak z_{\mathfrak p}(\mathfrak b)$ is the centralizer of $\mathfrak b$ 
in $\mathfrak p$.}

\vspace{0.5truecm}

\noindent
{\it Remark 1.2.} Since 
$\sharp(\overline{\triangle}_+\setminus\overline{\triangle}_+^1)\times 2+
\sharp\overline{\triangle}_+^1\leq \sharp(\triangle_+\setminus\triangle_+^1)
\times 2+\sharp\triangle_+^1$ (where 
$\triangle_+^1:=\{\alpha\in\triangle_+\,\vert\,{\rm the}\,\,{\rm multiplicity}
\,\,{\rm of}\,\,\alpha\,\,{\rm is}\,\,{\rm equal}\,\,{\rm to}\,\,1\}$), 
we have 
$$\sharp{\rm Spec}\,A_v\leq\sharp(\triangle_+\setminus\triangle_+^1)\times 2
+\sharp\triangle_+^1+{\rm dim}\,{\mathfrak z}_{\mathfrak p}({\mathfrak b})
-{\rm codim}\,M.  \leqno{(1.1)}$$

\vspace{0.5truecm}

In particular, we have the following fact.  

\vspace{0.3truecm}

\noindent
{\bf Corollary C.1.} {\sl Let $M$ be as in Theorem C.  Assume that 
${\rm codim}\,M={\rm rank}(G/K)$.  Then, for each normal vector $v$ of $M$, 
we have $\sharp{\rm Spec}\,A_v\leq\sharp(\triangle_+\setminus\triangle^1_+)
\times 2+\sharp\triangle^1_+$, where 
$\triangle_+^1:=\{\alpha\in\triangle_+\,\vert\,
{\rm the}\,\,{\rm multiplicity}\,\,{\rm of}\,\,\alpha\,\,{\rm is}\,\,
{\rm equal}\,\,{\rm to}\,\,1\}$.}

\vspace{0.3truecm}

In Table 1, we list up the number 
$m_{G/K}:=\sharp(\triangle_+\setminus\triangle^1_+)\times 2
+\sharp\triangle^1_+$ for irreducible symmetric spaces $G/K$'s 
of non-compact type.  Also, in Appendix 1, we list up the numbers 
${\rm max}_{v\in T^{\perp}M}\sharp{\rm Spec}\,A_v$ for principal orbits of 
Hermann type actions $H$'s on irreducible symmetric spaces $G/K$'s of 
non-compact type satisfying ${\rm cohom}\,H={\rm rank}(G/K)$.  

\vspace{0.3truecm}

\noindent
{\bf Future plan of research.} 
{\sl By using Theorems B and C, we will investigate whether the above question 
is solved affirmatively in some symmetric spaces of non-compact type.}

\vspace{0.3truecm}

For the focal set of a curvature-adapted and proper complex equifocal 
$C^{\omega}$-submanifold, we have the following fact.  

\vspace{0.3truecm}

\noindent
{\bf Theorem D.} {\sl Let $M$ be as in Theorem B.  Then the focal set of 
$(M,x_0)$ ($x_0\,:\,$ an arbitrary point of $M$) consists of finitely many 
totally geodesic hypersurfaces through some point in the section 
$\Sigma:=\exp^{\perp}(T^{\perp}_{x_0}M)$.}

\vspace{0.3truecm}

Let $\{{\it l}_i\,\vert\,i=1,\cdots,k\}$ be hyperplanes of 
$T^{\perp}_{x_0}M$ such that $\displaystyle{\mathop{\cup}_{i=1}^k\exp^{\perp}
({\it l}_i)}$ is the focal set of $(M,x_0)$.  Denote by $W_{M,{\bf R}}$ 
the group generated by the reflections with respect to ${\it l}_i$'s 
($i=1,\cdots,k$).  In this paper, we call this group the 
{\it real Coxeter group associated with} $M$ ({\it at} $x_0$).  Note that 
this group is independent of the choice of the base point $x_0$ up to 
isomorphicness.  For this group, we have the following fact.  

\vspace{0.3truecm}

\noindent
{\bf Theorem E.} {\sl Let $M$ and 
$\overline{\triangle}$ be as in Theorem B.  Then 
the real Coxeter group associated with $M$ is isomorphic to a subgroup of 
the Coxeter group $W_{\overline{\triangle}}$.}

\vspace{0.3truecm}

\noindent
{\it Remark 1.3.} We consider the case where $M$ is a principal orbit of 
a Hermann type action $H\times G/K\to G/K$.  Let $\sigma$ (resp. $\theta$) 
be an involution of $G$ with $({\rm Fix}\,\sigma)_0\subset H\subset{\rm Fix}\,
\sigma$ (resp. $({\rm Fix}\,\theta)_0\subset K\subset{\rm Fix}\,\theta$), 
where we may assume $\sigma\circ\theta=\theta\circ\sigma$ without loss of 
generality.  Then the real Coxeter group associated with $M$ is isomorphic to 
the Weyl group associated with the symmetric space 
${\rm Fix}(\sigma\circ\theta)/H\cap K$ (see Appendix 1).  

\section{Basic notions and facts}
In this section, we recall basic notions introduced in [Koi1$\sim$4].  
We first recall the notion of a complex equifocal submanifold.  
Let $M$ be an immersed submanifold with abelian normal bundle 
in a symmetric space $N=G/K$ of non-compact type.  Denote by $A$ the shape 
tensor of $M$.  Let $v\in T^{\perp}_xM$ and $X\in T_xM$ ($x=gK$).  Denote 
by $\gamma_v$ the geodesic in $N$ with $\dot{\gamma}_v(0)=v$.  
The strongly $M$-Jacobi field $Y$ along $\gamma_v$ with $Y(0)=X$ (hence 
$Y'(0)=-A_vX$) is given by 
$$Y(s)=(P_{\gamma_v\vert_{[0,s]}}\circ(D^{co}_{sv}-sD^{si}_{sv}\circ A_v))
(X),$$

\newpage

$$\begin{tabular}{|l|c|c|c|c|}
\hline
{\rm Type} & $G/K$ & $\sharp\triangle_+$ & $\sharp\triangle^1_+$ & 
$m_{G/K}$\\
\hline
{\scriptsize{\rm(AI)}} & {\scriptsize$SL(n,{\bf R})/SO(n)\,\,(n\geq3)$} & 
{\scriptsize$\frac{n(n-1)}{2}$} & {\scriptsize$\frac{n(n-1)}{2}$} & 
{\scriptsize$\frac{n(n-1)}{2}$}\\
\hline
{\scriptsize{\rm(AII)}} & {\scriptsize$SU^{\ast}(2n)/Sp(n)\,\,(n\geq3)$} & 
{\scriptsize$\frac{n(n-1)}{2}$} & {\scriptsize$0$} & {\scriptsize$n(n-1)$}\\
\hline
{\scriptsize{\rm(AIII)}} & {\scriptsize$\displaystyle{
\begin{array}{r}
SU(p,q)/S(U(p)\times U(q))\\
(1\leq p<q)
\end{array}}$}
& {\scriptsize$p^2+p$} & {\scriptsize$p$} & 
{\scriptsize$p(2p+1)$}\\
 & {\scriptsize$\displaystyle{
\begin{array}{r}
SU(p,p)/S(U(p)\times U(p))\\
(p\geq2)
\end{array}}$}
& {\scriptsize$p^2$} & {\scriptsize$p$} & {\scriptsize$p(2p-1)$}\\
\hline
{\scriptsize{\rm(BDI)}} & 
{\scriptsize$\displaystyle{
\begin{array}{r}
SO_0(p,q)/SO(p)\times SO(q)\\
(2\leq p<q)
\end{array}}$} & {\scriptsize$p^2$} & 
{\scriptsize$\displaystyle{
\left\{\begin{array}{c}
p^2\\
p(p-1)
\end{array}
\right.}$}
&
{\scriptsize$\displaystyle{
\begin{array}{cl}
p^2 & (q-p=1)\\
p(p+1) & (q-p\geq2)
\end{array}
}$}\\
{\scriptsize } & 
{\scriptsize$
SO_0(1,q)/SO(1)\times SO(q)$} & {\scriptsize$1$} & 
{\scriptsize$\displaystyle{
\left\{\begin{array}{c}
1\\
0
\end{array}
\right.}$}
&
{\scriptsize$\displaystyle{
\begin{array}{cl}
1 & (q=2)\\
2 & (q\geq3)
\end{array}
}$}\\
\hline
{\scriptsize{\rm(BDI$'$)}} & {\scriptsize$SO_0(p,p)/SO(p)\times SO(p)$} & 
{\scriptsize$p(p-1)$} & {\scriptsize$p(p-1)$} & {\scriptsize$p(p-1)$}\\
\hline
{\scriptsize{\rm(DIII)}} & 
{\scriptsize
$\displaystyle{
\begin{array}{r}
SO^{\ast}(2n)/U(n)\\
(n\geq4)
\end{array}}$} & 
{\scriptsize$\displaystyle{
\left\{\begin{array}{c}
\frac{n^2-1}{4}\\
\frac{n^2}{4}
\end{array}
\right.}$}
& 
{\scriptsize$\displaystyle{
\begin{array}{c}
\frac{n-1}{2}\\
\frac{n}{2}
\end{array}}$}
&
{\scriptsize$\displaystyle{
\begin{array}{cl}
\frac{n(n-1)}{2} & (n:{\rm odd})\\
\frac{n(n-1)}{2} & (n:{\rm even})
\end{array}
}$}\\
\hline
{\scriptsize{\rm(CI)}} & {\scriptsize$Sp(n,{\bf R})/U(n)\,\,(n\geq2)$} & 
{\scriptsize$n^2$} & 
{\scriptsize$n^2$} & {\scriptsize$n^2$}\\
\hline
{\scriptsize{\rm(CII)}} & 
{\scriptsize$\displaystyle{
\begin{array}{r}
Sp(p,q)/Sp(p)\times Sp(q)\\
(p<q)
\end{array}}$} & {\scriptsize$p(p+1)$} & {\scriptsize$0$} & 
{\scriptsize$2p(p+1)$}\\
& {\scriptsize$\displaystyle{
\begin{array}{r}
Sp(p,p)/Sp(p)\times Sp(p)\\
(p\geq2)
\end{array}}$} & {\scriptsize$p^2$} & {\scriptsize$0$} & 
{\scriptsize$2p^2$}\\
\hline
{\scriptsize{\rm(EI)}} & {\scriptsize$E_6^6/Sp(4)$} & {\scriptsize$36$} & 
{\scriptsize$36$} & {\scriptsize$36$}\\
\hline
{\scriptsize{\rm(EII)}} & {\scriptsize$E_6^2/SU(6)\cdot SU(2)$} & {\scriptsize$24$} & 
{\scriptsize$12$} & {\scriptsize$36$}\\
\hline
{\scriptsize{\rm(EIII)}} & {\scriptsize$E_6^{-14}/Spin(10)\cdot U(1)$} & 
{\scriptsize$6$} & {\scriptsize$2$} & {\scriptsize$10$}\\
\hline
{\scriptsize{\rm(EIV)}} & {\scriptsize$E_6^{-26}/F_4$} & {\scriptsize$3$} & 
{\scriptsize$0$} & {\scriptsize$6$}\\
\hline
{\scriptsize{\rm(EV)}} & {\scriptsize$E_7^7/(SU(8)/\{\pm1\})$} & 
{\scriptsize$63$} & {\scriptsize$63$} & {\scriptsize$63$}\\
\hline
{\scriptsize{\rm(EVI)}} & {\scriptsize$E_7^{-5}/SO'(12)\cdot SU(2)$} & 
{\scriptsize$24$} 
& {\scriptsize$12$} & {\scriptsize$36$}\\
\hline
{\scriptsize{\rm(EVII)}} & {\scriptsize$E_7^{-25}/E_6\cdot U(1)$} & 
{\scriptsize$9$} 
& {\scriptsize$3$} & {\scriptsize$15$}\\
\hline
{\scriptsize{\rm(EVIII)}} & {\scriptsize$E_8^8/SO'(16)$} & {\scriptsize$120$} 
& {\scriptsize$120$} & {\scriptsize$120$}\\
\hline
{\scriptsize{\rm(EIX)}} & {\scriptsize$E_8^{-24}/E_7\cdot Sp(1)$} & 
{\scriptsize$24$} & 
{\scriptsize$12$} & {\scriptsize$36$}\\
\hline
{\scriptsize{\rm(FI)}} & {\scriptsize$F_4^4/Sp(3)\cdot Sp(1)$} & 
{\scriptsize$24$} & {\scriptsize$24$} & {\scriptsize$24$}\\
\hline
{\scriptsize{\rm(FII)}} & {\scriptsize$F_4^{-20}/Spin(9)$} & {\scriptsize$2$} 
& {\scriptsize$0$} & {\scriptsize$4$}\\
\hline
{\scriptsize{\rm(G)}} & {\scriptsize$G_2^2/SO(4)$} & {\scriptsize$6$} 
& {\scriptsize$6$} & {\scriptsize$6$}\\
\hline
{\scriptsize{\rm(II-A)}} & {\scriptsize$\displaystyle{
\begin{array}{r}
SL(n,{\bf C})/SU(n)\\
(n\geq3)
\end{array}}$}
& {\scriptsize$\frac{n(n-1)}{2}$} & {\scriptsize$0$} & 
{\scriptsize$n(n-1)$}\\
\hline
{\scriptsize{\rm(II-BD)}} & {\scriptsize$\displaystyle{
\begin{array}{r}
SO(n,{\bf C})/SO(n)\\
(n\geq4)
\end{array}}$}

& 
{\scriptsize$\displaystyle{
\left\{\begin{array}{c}
\frac{(n-1)^2}{4}\\
\frac{n(n-2)}{4}
\end{array}\right.}$} & 
{\scriptsize$\displaystyle{
\begin{array}{c}
0\\
0
\end{array}}$} & 
{\scriptsize$\displaystyle{
\begin{array}{cl}
\frac{(n-1)^2}{2} & (n:{\rm odd})\\
\frac{n(n-2)}{2} & (n:{\rm even})
\end{array}}$}\\
\hline
{\scriptsize{\rm(II-C)}} & {\scriptsize$Sp(n,{\bf C})/Sp(n)$} & 
{\scriptsize$n^2$} & {\scriptsize$0$} & 
{\scriptsize$2n^2$}\\
\hline
{\scriptsize{\rm(II-}${\rm E}_6$)} & {\scriptsize$E_6^{\bf c}/E_6$} 
& {\scriptsize$36$} & 
{\scriptsize$0$} & {\scriptsize$72$}\\
\hline
{\scriptsize{\rm(II-}${\rm E}_7$)} & {\scriptsize$E_7^{\bf c}/E_7$} 
& {\scriptsize$63$} & 
{\scriptsize$0$} & {\scriptsize$126$}\\
\hline
{\scriptsize{\rm(II-}${\rm E}_8$)} & {\scriptsize$E_8^{\bf c}/E_8$} 
& {\scriptsize$120$} & 
{\scriptsize$0$} & {\scriptsize$240$}\\
\hline
{\scriptsize{\rm(II-}${\rm F}_4$)} & {\scriptsize$F_4^{\bf c}/F_4$} 
& {\scriptsize$24$} & 
{\scriptsize$0$} & {\scriptsize$48$}\\
\hline
{\scriptsize{\rm(II-}${\rm G}_2$)} & {\scriptsize$G_2^{\bf c}/G_2$} 
& {\scriptsize$6$} & 
{\scriptsize$0$} & {\scriptsize$12$}\\
\hline
\end{tabular}$$

\vspace{0.2truecm}

\centerline{{\bf Table 1.}}

\newpage

\noindent
where $Y'(0)=\widetilde{\nabla}_vY,\,\,P_{\gamma_v\vert_{[0,s]}}$ is 
the parallel translation along $\gamma_v\vert_{[0,s]}$ and 
$D^{co}_{sv}$ (resp. $D^{si}_{sv}$) is given by 
$$\begin{array}{c}
\displaystyle{
D^{co}_{sv}=g_{\ast}\circ\cos(\sqrt{-1}{\rm ad}(sg_{\ast}^{-1}v))
\circ g_{\ast}^{-1}}\\
\displaystyle{\left({\rm resp.}\,\,\,\,
D^{si}_{sv}=g_{\ast}\circ
\frac{\sin(\sqrt{-1}{\rm ad}(sg_{\ast}^{-1}v))}
{\sqrt{-1}{\rm ad}(sg_{\ast}^{-1}v)}\circ g_{\ast}^{-1}\right).}
\end{array}$$ 
Here ${\rm ad}$ is the adjoint 
representation of the Lie algebra $\mathfrak g$ of $G$.  
All focal radii of $M$ along $\gamma_v$ are obtained as real numbers $s_0$ 
with ${\rm Ker}(D^{co}_{s_0v}-s_0D^{si}_{s_0v}\circ A_v)\not=\{0\}$.  So, we 
call a complex number $z_0$ with ${\rm Ker}(D^{co}_{z_0v}-
z_0D^{si}_{z_0v}\circ A_v^{{\bf c}})\not=\{0\}$ a {\it complex 
focal radius of} $M$ {\it along} $\gamma_v$ and call ${\rm dim}\,
{\rm Ker}(D^{co}_{z_0v}-z_0D^{si}_{z_0v}\circ A_v^{{\bf c}})$ the 
{\it multiplicity} of the complex focal radius $z_0$, 
where $A_v^{\bf c}$ is the complexification of $A_v$ and $D^{co}_{z_0v}$ 
(resp. $D^{si}_{z_0v}$) is a ${\bf C}$-linear transformation of 
$(T_xN)^{\bf c}$ defined by 
$$\begin{array}{c}
\displaystyle{
D^{co}_{z_0v}=g^{\bf c}_{\ast}\circ\cos(\sqrt{-1}{\rm ad}^{\bf c}
(z_0g_{\ast}^{-1}v))\circ (g^{\bf c}_{\ast})^{-1}}\\
\displaystyle{\left({\rm resp.}\,\,\,\,
D^{si}_{sv}=g^{\bf c}_{\ast}\circ
\frac{\sin(\sqrt{-1}{\rm ad}^{\bf c}(z_0g_{\ast}^{-1}v))}
{\sqrt{-1}{\rm ad}^{\bf c}(z_0g_{\ast}^{-1}v)}\circ(g^{\bf c}_{\ast})^{-1}
\right),}
\end{array}$$
where $g_{\ast}^{\bf c}$ (resp. ${\rm ad}^{\bf c}$) is the complexification 
of $g_{\ast}$ (resp. ${\rm ad}$).  
Here we note that, in the case where $M$ is of class $C^{\omega}$, 
complex focal radii along $\gamma_v$ 
indicate the positions of focal points of the extrinsic 
complexification $M^{\bf c}(\hookrightarrow G^{\bf c}/K^{\bf c})$ of $M$ 
along the complexified geodesic $\gamma_{\iota_{\ast}v}^{\bf c}$, where 
$G^{\bf c}/K^{\bf c}$ is the anti-Kaehler symmetric space associated with 
$G/K$ and $\iota$ is the natural immersion of $G/K$ into 
$G^{\bf c}/K^{\bf c}$.  
See Section 4 of [Koi2] about the definitions of 
$G^{\bf c}/K^{\bf c},\,
M^{\bf c}(\hookrightarrow G^{\bf c}/K^{\bf c})$ and 
$\gamma_{\iota_{\ast}v}^{\bf c}$.  
Also, for a complex focal radius $z_0$ of $M$ along $\gamma_v$, we 
call $z_0v$ ($\in (T^{\perp}_xM)^{\bf c}$) a 
{\it complex focal normal vector of} $M$ {\it at} $x$.  
Furthermore, assume that $M$ has globally flat normal bundle, that is, 
the normal holonomy group of $M$ is trivial.  
Let $\tilde v$ be a parallel unit normal vector field of $M$.  
Assume that the number (which may be $0$ and $\infty$) of distinct complex 
focal radii along $\gamma_{\tilde v_x}$ is independent of the choice of 
$x\in M$.  Furthermore assume that the number is not equal to $0$.  
Let $\{r_{i,x}\,\vert\,i=1,2,\cdots\}$ 
be the set of all complex focal radii along $\gamma_{\tilde v_x}$, where 
$\vert r_{i,x}\vert\,<\,\vert r_{i+1,x}\vert$ or 
"$\vert r_{i,x}\vert=\vert r_{i+1,x}\vert\,\,\&\,\,{\rm Re}\,r_{i,x}
>{\rm Re}\,r_{i+1,x}$" or 
"$\vert r_{i,x}\vert=\vert r_{i+1,x}\vert\,\,\&\,\,
{\rm Re}\,r_{i,x}={\rm Re}\,r_{i+1,x}\,\,\&\,\,
{\rm Im}\,r_{i,x}=-{\rm Im}\,r_{i+1,x}<0$".  
Let $r_i$ ($i=1,2,\cdots$) be complex 
valued functions on $M$ defined by assigning $r_{i,x}$ to each $x\in M$.  
We call these functions $r_i$ ($i=1,2,\cdots$) {\it complex 
focal radius functions for} $\tilde v$.  
We call $r_i\tilde v$ a {\it complex focal normal vector field for} 
$\tilde v$.  If, for each parallel 
unit normal vector field $\tilde v$ of $M$, the number of distinct complex 
focal radii along $\gamma_{\tilde v_x}$ is independent of the choice of 
$x\in M$, each complex focal radius function for $\tilde v$ 
is constant on $M$ and it has constant multiplicity, then 
we call $M$ a {\it complex equifocal submanifold}.  

Next we shall recall the notion of a proper complex equifocal submanifold.  
For its purpose, we first recall the notion of a proper complex isoparametric 
submanifold in a pseudo-Hilbert space.  
Let $M$ be a pseudo-Riemannian Hilbert submanifold in a pseudo-Hilbert space 
$(V,\langle\,\,,\,\,\rangle)$ immersed by $f$.  
See Section 2 of [Koi1] about the definitions of a pseudo-Hilbert space and 
a pseudo-Riemannian Hilbert submanifold.  
Denote by $A$ the shape tensor of $M$ and by $T^{\perp}M$ the normal bundle 
of $M$.  Note that, for $v\in T^{\perp}M$, $A_v$ is not necessarily 
diagonalizable with respect to an orthonormal base.  
We  call $M$ a {\it Fredholm pseudo-Riemannian Hilbert submanifold} 
(or simply {\it Fredholm submanifold}) if the following conditions hold:

\vspace{0.2truecm}

(F-i) $M$ is of finite codimension,

(F-ii) There exists an orthogonal time-space decomposition 
$V=V_-\oplus V_+$ such that $(V,\langle\,\,,\,\,\rangle_{V_{\pm}})$ is 
a Hilbert space and that, for each $v\in T^{\perp}M$, $A_v$ is a compact 
operator with respect to $f^{\ast}\langle\,\,,\,\,\rangle_{V_{\pm}}$.  

\vspace{0.2truecm}

\noindent
Since $A_v$ is a compact operator with respect to 
$f^{\ast}\langle\,\,,\,\,\rangle_{V_{\pm}}$, 
the operator ${\rm id}-A_v$ is a Fredholm operator with respect to 
$f^{\ast}\langle\,\,,\,\,\rangle_{V_{\pm}}$ and hence the normal exponential 
map $\exp^{\perp}\,:\,T^{\perp}M\to V$ of 
$M$ is a Fredholm map with respect to the metric of $T^{\perp}M$ naturally 
defined from $f^{\ast}\langle\,\,,\,\,\rangle_{V_{\pm}}$ and 
$\langle\,\,,\,\,\rangle_{V_{\pm}}$, where ${\rm id}$ is the identity 
transformation of $TM$.  
The spectrum of the complexification 
$A_v^{\bf c}$ of $A_v$ is described as 
$\{0\}\cup\{\lambda_i\,\vert\,i=1,2,\cdots\}$, where "$\vert\lambda_i\vert\,>\,
\vert\lambda_{i+1}\vert$" or "$\vert\lambda_i\vert=\vert\lambda_{i+1}\vert\,\,
\&\,\,{\rm Re}\,\lambda_i>{\rm Re}\,\lambda_{i+1}$" or 
"$\vert\lambda_i\vert=\vert\lambda_{i+1}\vert\,\,\&\,\,
{\rm Re}\,\lambda_i={\rm Re}\,\lambda_{i+1}\,\,\&\,\,{\rm Im}\,\lambda_i
=-{\rm Im}\,\lambda_{i+1}>0$".  We call 
$\lambda_i$ the $i$-{\it th} {\it complex principal curvature of direction} 
$v$.  Assume that $M$ has globally flat normal bundle.  
Fix a parallel normal vector field $\tilde v$ on $M$.  
Assume that the number (which may be $\infty$) 
of distinct complex principal curvatures of $\tilde v_x$ is 
independent of the choice of $x\in M$.  Then we can define functions 
$\tilde{\lambda}_i$ 
($i=1,2,\cdots$) on $M$ by assigning the $i$-th complex principal curvature 
of direction $\tilde v_x$ to each $x\in M$.  We call this function 
$\tilde{\lambda}_i$ the $i$-{\it th} {\it complex principal curvature function 
of direction} $\tilde v$.  
If $M$ is a Fredholm submanifold with globally flat normal bundle 
satisfying the following condition (CI), then we call $M$ 
a {\it complex isoparametric submanifold}:

\vspace{0.2truecm}

\noindent
(CI) $\quad$ for each parallel normal vector field $\tilde v$, 
the number of distinct complex principal curvatures of direction 
$\tilde v_x$ is independent of the choice of $x\in M$ and 
each complex principal curvature function of direction $\tilde v$ 
is constant on $M$ and has constant multiplicity.  

\vspace{0.2truecm}

\noindent
Furthermore, if, for each $v\in T^{\perp}M$, there exists a pseudo-orthonormal 
base of $(T_xM)^{\bf c}$ ($x\,:\,$ the base point of $v$) consisting 
of the eigenvectors of the complexified shape operator $A_v^{\bf c}$, 
then we  call $M$ a {\it proper complex isoparametric submanifold}.  
Then, for each $x\in M$, there exists a pseudo-orthonormal base 
of $(T_xM)^{\bf c}$ consisting of the common-eigenvectors of 
the complexified shape operators 
$A_v^{\bf c}$'s ($v\in T^{\perp}_xM$) because 
$A_v^{\bf c}$'s are commutative.  
Let $\{E_i\,\vert\,i\in I\}$ ($I\subset {\bf N}$) be the family of 
subbundles of $(TM)^{\bf c}$ such that, for each $x\in M$, 
$\{E_i(x)\,\vert\,i\in I\}$ is the set of all common-eigenspaces 
of $A_v^{\bf c}$'s ($v\in T_x^{\perp}M$).  
Note that $\displaystyle{(T_xM)^{\bf c}=
\overline{\mathop{\oplus}_{i\in I}E_i(x)}}$ holds.  
There exist smooth sections $\lambda_i$ ($i\in I$) of 
$((T^{\perp}M)^{\bf c})^{\ast}$ such that $A_v^{\bf c}=\lambda_i(v){\rm id}$ 
on $E_i(\pi(v))$ for each $v\in T^{\perp}M$, where 
$\pi$ is the bundle projection of $(T^{\perp}M)^{\bf c}$.  We call 
$\lambda_i$ ($i\in I$) {\it complex principal curvatures of} $M$ and call 
subbundles $E_i$ ($i\in I$) of $(T^{\perp}M)^{\bf c}$ 
{\it complex curvature distributions of} $M$.  Note that 
$\lambda_i(v)$ is one of the complex principal curvatures of direction $v$.  
Set ${\it l}_i:=\lambda_i^{-1}(1)\,(\subset(T^{\perp}_xM)^{\bf c})$ and $R_i$ 
be the complex reflection of order two with respect to ${\it l}_i$, where 
$i\in I$.  Denote by $W_M$ the group generated by $R_i$'s ($i\in I$) which is 
independent of the choice of $x\in M$ up to isomorphicness.  We call 
${\it l}_i$'s {\it complex focal hyperplanes of} $(M,x)$.  
Let $N=G/K$ be a symmetric space of non-compact type and $\pi$ be the 
natural projection of $G$ onto $G/K$.  
Let $(\frak g,\sigma)$ be the orthogonal symmetric Lie algebra of $G/K$, 
$\mathfrak f=\{X\in \mathfrak g\,\vert\,\sigma(X)=X\}$ and 
$\mathfrak p=\{X\in \mathfrak g\,\vert\,\sigma(X)=-X\}$, which is identified 
with the tangent space $T_{eK}N$.  Let $\langle\,\,,\,\,\rangle$ be 
the ${\rm Ad}(G)$-invariant non-degenerate symmetric bilinear form of 
$\mathfrak g$ inducing the Riemannian metric of $N$.  Note that 
$\langle\,\,,\,\,\rangle\vert_{\mathfrak f\times\mathfrak f}$ 
(resp. $\langle\,\,,\,\,\rangle\vert_{\mathfrak p\times\mathfrak p}$) is 
negative (resp. positive) definite.  Denote by the same symbol 
$\langle\,\,,\,\,\rangle$ the bi-invariant pseudo-Riemannian metric of $G$ 
induced from $\langle\,\,,\,\,\rangle$ and the Riemannian metric of $N$.  
Set $\mathfrak g_+:=\mathfrak p,\,\,\mathfrak g_-:=\mathfrak f$ and 
$\langle\,\,,\,\,\rangle_{\mathfrak g_{\pm}}
:=-\pi_{\mathfrak g_-}^{\ast}\langle\,\,,\,\,\rangle
+\pi_{\mathfrak g_+}^{\ast}\langle\,\,,\,\,\rangle$, where 
$\pi_{\mathfrak g_-}$ (resp. $\pi_{\mathfrak g_+}$) is the projection of 
$\mathfrak g$ onto $\mathfrak g_-$ (resp. $\mathfrak g_+$).  Let 
$H^0([0,1],\mathfrak g)$ be the space of all $L^2$-integrable paths 
$u:[0,1]\to\frak g$ (with respect to 
$\langle\,\,,\,\,\rangle_{\mathfrak g_{\pm}}$).  
Let $H^0([0,1],\mathfrak g_-)$ (resp. $H^0([0,1],\mathfrak g_+)$) be 
the space of all $L^2$-integrable paths 
$u:[0,1]\to\mathfrak g_-$ (resp. $u:[0,1]\to \mathfrak g_+$) with respect to 
$-\langle\,\,,\,\,\rangle\vert_{\mathfrak g_-\times\mathfrak g_-}$ 
(resp. $\langle\,\,,\,\,\rangle\vert_{\mathfrak g_+\times\mathfrak g_+}$).  
It is clear that $H^0([0,1],\mathfrak g)=H^0([0,1],\mathfrak g_-)\oplus 
H^0([0,1],\mathfrak g_+)$.  
Define a non-degenerate symmetric bilinear form 
$\langle\,\,,\,\,\rangle_0$ of 
$H^0([0,1],\mathfrak g)$ by 
$\langle u,v\rangle_0:=\int_0^1\langle u(t),v(t)\rangle dt$.  It is easy to 
show that the decomposition $H^0([0,1],\mathfrak g)=H^0([0,1],\mathfrak g_-)
\oplus H^0([0,1],\mathfrak g_+)$ is an orthogonal time-space decomposition 
with respect to $\langle\,\,,\,\,\rangle_0$.  For simplicity, set 
$H^0_{\pm}:=H^0([0,1],\mathfrak g_{\pm})$ and 
$\langle\,\,,\,\,\rangle_{0,H^0_{\pm}}:=-\pi^{\ast}_{H^0_-}
\langle\,\,,\,\,\rangle_0+\pi^{\ast}_{H^0_+}\langle\,\,,\,\,\rangle_0$, where 
$\pi_{H^0_-}$ (resp. $\pi_{H^0_+}$) is the projection of 
$H^0([0,1],\mathfrak g)$ onto $H^0_-$ (resp. $H^0_+$).  
It is clear that $\langle u,v\rangle_{0,H^0_{\pm}}
=\int_0^1\langle u(t),v(t)\rangle_{\mathfrak g_{\pm}}dt$ 
($u,\,v\in H^0([0,1],\mathfrak g)$).  
Hence $(H^0([0,1],\mathfrak g),\,\langle\,\,,\,\,\rangle_{0,H^0_{\pm}})$ is a 
Hilbert space, that is, $(H^0([0,1],\mathfrak g),\,\langle\,\,,\,\,\rangle_0)$ 
is a pseudo-Hilbert space.  Let $H^1([0,1],G)$ be the Hilbert Lie group of all 
absolutely continuous paths $g:[0,1]\to G$ such that the weak derivative $g'$ 
of $g$ is squared integrable (with respect to 
$\langle\,\,,\,\,\rangle_{\mathfrak g_{\pm}}$), that is, 
$g_{\ast}^{-1}g'\in H^0([0,1],\mathfrak g)$.  Define a map 
$\phi:H^0([0,1],\mathfrak g)\to G$ by $\phi(u)=g_u(1)$ 
($u\in H^0([0,1],\mathfrak g)$), where $g_u$ is the element of 
$H^1([0,1],G)$ satisfying $g_u(0)=e$ and $g_{u\ast}^{-1}g_u'=u$.  
We call this map the {\it parallel transport map} (from $0$ to $1$).  
This submersion $\phi$ is a pseudo-Riemannian submersion of 
$(H^0([0,1],\mathfrak g),\langle\,\,,\,\,\rangle_0)$ onto 
$(G,\langle\,\,,\,\,\rangle)$.  
Let $\pi:G\to G/K$ be the natural projection.  It follows from Theorem A of 
[Koi1] (resp. Theorem 1 of [Koi2]) that, in the case where $M$ is curvature 
adapted (resp. of class $C^{\omega}$), $M$ is complex equifocal if and only if 
each component of $(\pi\circ\phi)^{-1}(M)$ is complex isoparametric.  
In particular, if components of $(\pi\circ\phi)^{-1}(M)$ 
are proper complex isoparametric, then we call $M$ 
a {\it proper complex equifocal submanifold}.  
Let $M$ be a proper complex equifocal $C^{\omega}$-submanifold in $G/K$, 
$\widetilde M_0$ be a component of $\widetilde M:=(\pi\circ\phi)^{-1}(M)$.  
Denote by $W_{\widetilde M_0}$ the group defined as above for this proper 
complex isoparametric submanifold $\widetilde M_0$, where we take $u_0$ as 
the base point.  

Let $N=G/K$ be a symmetric space of non-compact type, 
$({\mathfrak g},\sigma)$ be the orthogonal symmetric Lie algebra 
associated with a symmetric pair $(G,K)$ and 
${\mathfrak g}={\mathfrak f}+{\mathfrak p}$ be the Cartan decomposition.  
Note that $\mathfrak f$ is the Lie algebra of $K$ and $\mathfrak p$ is 
identified with the tangent space $T_{eK}N$, where $e$ is the identity 
element of $G$.  Let $\langle\,\,,\,\,\rangle$ be the 
${\rm Ad}(G)$-invariant non-degenerate inner product of $\mathfrak g$ 
inducing the Riemannian metric of $N$.  Let ${\mathfrak g}^{{\bf c}},
\,\,{\mathfrak f}^{{\bf c}},\,\,{\mathfrak p}^{{\bf c}}$ and 
$\langle\,\,,\,\,\rangle^{{\bf c}}$ be the complexifications of 
$\mathfrak g,\,\,\mathfrak f,\,\,\mathfrak p$ and $\langle\,\,,\,\,\rangle$, 
respectively.  
Let $\mathfrak a$ be a maximal abelian subspace of $\mathfrak p$ and 
$\mathfrak p=\mathfrak a+\sum\limits_{\alpha\in\triangle_+}
\mathfrak p_{\alpha}$ be 
the root space decomposition with respect to $\mathfrak a$.  
Then $(\mathfrak g^{\bf c},\mathfrak f^{\bf c})$ is a semi-simple symmetric 
pair, $\mathfrak a$ is a maximal split abelian subspace 
of $\mathfrak p^{\bf c}$ and 
$\mathfrak p^{{\bf c}}=\mathfrak a^{{\bf c}}+
\sum\limits_{\alpha\in\triangle_+}\mathfrak p_{\alpha}^{{\bf c}}$ is the root 
space decomposition with respect to $\mathfrak a$, where 
$\mathfrak a^{{\bf c}}$ and $\mathfrak p_{\alpha}^{{\bf c}}$ are 
the complexifications of $\mathfrak a$ and $\mathfrak p_{\alpha}$, 
respectively.  
Note that $\mathfrak a^{\bf c}$ is the centralizer of $\mathfrak a$ in 
$\mathfrak p^{\bf c}$.  See [R] and [OS] about the general theory of 
a semi-simple symmetric pair.  
Let $G^{{\bf c}}$ (resp. $K^{{\bf c}}$) be the complexification 
of $G$ (resp. $K$).  The $2$-multiple of the real part 
${\rm Re}\langle\,\,,\,\,\rangle^{{\bf c}}$ of 
$\langle\,\,,\,\,\rangle^{{\bf c}}$ is the Killing form of 
$\mathfrak g^{\bf c}$ regarded as a real Lie algebra.  The restriction 
$2{\rm Re}\langle\,\,,\,\,\rangle^{\bf c}\vert_{{\mathfrak p}^{\bf c}\times
{\mathfrak p}^{\bf c}}$ is an ${\rm Ad}(K^{\bf c})$-
invariant non-degenerate inner product of ${\mathfrak p}^{\bf c}$ 
($=T_{eK^{\bf c}}(G^{\bf c}/K^{\bf c})$).  
Denote by $\langle\,\,,\,\,\rangle'$ the $G^{{\bf c}}$-invariant 
pseudo-Riemannian metric on $G^{{\bf c}}/K^{{\bf c}}$ 
induced from $2{\rm Re}\langle\,\,,\,\,\rangle^{{\bf c}}
\vert_{{\mathfrak p}^{{\bf c}}\times{\mathfrak p}^{{\bf c}}}$.  
Define an almost complex structure $J_0$ of ${\mathfrak p}^{{\bf c}}$ 
by $J_0(X+\sqrt{-1}Y)=-Y+\sqrt{-1}X$ ($X,Y\in\mathfrak p$).  It is clear that 
$J_0$ is ${\rm Ad}(K^{{\bf c}})$-invariant.  Denote by 
$\widetilde J$ the $G^{{\bf c}}$-invariant almost complex structure 
on $G^{{\bf c}}/K^{{\bf c}}$ induced from $J_0$.  It is shown 
that $(G^{{\bf c}}/K^{{\bf c}},\langle\,\,,\,\,\rangle',
\widetilde J)$ is an anti-Kaehlerian manifold and a (semi-simple) 
pseudo-Riemannian symmetric space.  We call this anti-Kaehlerian 
manifold an {\it anti-Kaehlerian symmetric space associated with} $G/K$ 
and simply denote it by $G^{{\bf c}}/K^{{\bf c}}$.  
Let $\pi^{\bf c}:G^{\bf c}\to G^{\bf c}/K^{\bf c}$ be the natural projection 
and $\phi^{\bf c}:H^0([0,1],\mathfrak g^{\bf c})\to G^{\bf c}$ be the parallel 
transport map for $G^{\bf c}$.  This map $\phi^{\bf c}$ is defined in similar 
to $\phi$ (see Section 6 of [Koi2] in detail).  
Let $M$ be a complete $C^{\omega}$-submanifold in $G/K$ and $M^{\bf c}$ be 
the extrinsic complexification of $M$.  
Let $\widetilde{M^{\bf c}_0}$ be a 
component of $\widetilde{M^{\bf c}}:=(\pi^{\bf c}\circ\phi^{\bf c})^{-1}
(M^{\bf c})$.  In [Koi4], we defined the complex Coxeter group associated with 
$M$ as the group generated by complex reflections of order two with respect to 
complex focal hyperplanes constructing the focal set of 
$\widetilde{M^{\bf c}_0}$ at an arbitrary fixed point $u_1$.  Denote by $W_M$ 
this group, which is discrete (see Proposition 3.7 of [Koi4]).  Since 
the complex focal hyperplanes of $\widetilde M_0$ at $u_0$ coincides with 
those of $\widetilde{M^{\bf c}_0}$ at $u_1$ under some identification of 
$(T^{\perp}_{u_0}\widetilde M_0)^{\bf c}$ with $T^{\perp}_{u_1}
(\widetilde{M^{\bf c}_0})$, we see that $W_{\widetilde M_0}$ is isomorphic to 
$W_M$.  

At the end of this section, we recall the notions of the Weyl group and 
the affine Weyl group associated with a root system.  Let $\triangle$ be 
a subset of the dual space ${\Bbb E}^{\ast}$ of a Euclidean space ${\Bbb E}$ 
consisting of non-zero vectors.  We consider the following three conditions:

\vspace{0.3truecm}

(i) If $\alpha,\,\beta\in\triangle$, then $s_{\alpha}(\beta)\in\triangle$, 
where $s_{\alpha}$ is the reflection with respect to $\alpha^{-1}(0)$,

(ii) If $\alpha,\,\beta\in\triangle$, then $\frac{2\langle\alpha,\beta\rangle}
{\langle\alpha,\alpha\rangle}\in{\bf Z}$,

(iii) If $\alpha,\,a\alpha\in\triangle$ ($a\in{\bf R}$), then $a=\pm1$.

\vspace{0.3truecm}

\noindent
If $\triangle$ satisfies the condition (i), then we call it a 
{\it weakly root system}.  
Here we note that, if $\triangle$ satisfies the conditions (i) and (iii), then 
it is called a root system in [Ka].  
If $\triangle$ satisfies the conditions (i) and (ii), then it is called 
a {\it root system} (see [He]).  
Furthermore, if $\triangle$ satisfies the condition (iii), then it is said 
to be {\it reduced}.  For a weakly root system $\triangle$, we denote by 
$W_{\triangle}$ the (finite) group generated by the reflection's 
with respect to $\alpha^{-1}(0)$'s ($\alpha\in\triangle$) and 
denote by $W^A_{\triangle}$ the affine transformation group generated by 
the reflections with respect to $\alpha^{-1}(j)$'s 
($\alpha\in\triangle,\,j\in{\bf Z}$).  
We call $W_{\triangle}$ the {\it linear transformation group associated with} 
$\triangle$ and $W^A_{\triangle}$ the {\it affine transformation group 
associated with} $\triangle$.  If $W_{\triangle}$ is finite, then we call 
$W_{\triangle}$ the {\it Coxeter group associated with} $\triangle$ and, if 
$\triangle$ is a root system, then $W_{\triangle}$ is called the 
{\it Weyl group associated with} $\triangle$.  Also, if $W^A_{\triangle}$ is 
discrete, then we call $W^A_{\triangle}$ the {\it affine Weyl group associated 
with} $\triangle$.  






\section{Focal points on the ideal boundary}
In this section, we introduce the notion of a focal point on the ideal 
boundary $N(\infty)$ 
for a submanifold $M$ in a Hadamard manifold $N$.  
Denote by $\widetilde{\nabla}$ the Levi-Civita 
connection of $N$ and $A$ the shape tensor of $M$.  
Let $\gamma_v:[0,\infty)\to N$. be the normal geodesic 
of $M$ of direction $v\in T^{\perp}_xM$.  If there exists a $M$-Jacobi field 
(resp. strongly $M$-Jacobi field) $Y$ along $\gamma_v$ satisfying 
$\lim\limits_{t\to\infty}\frac{\vert\vert Y_t\vert\vert}{t}=0$, then we call 
$\gamma_v(\infty)\,(\in N(\infty))$ 
a {\it focal point} (resp. {\it strongly focal point}) 
{\it on the ideal boundary} $N(\infty)$ {\it of} $M$ {\it along} $\gamma_v$, 
where $\gamma_v(\infty)$ is the asymptotic class of $\gamma_v$ (see Fig. 1).  
Here a $M$-Jacobi field along $\gamma_v$ implies a 
Jacobi field $Y$ along $\gamma_v$ satisfying $Y(0)\in T_xM$ and $Y'(0)_T=
-A_vY(0)$ and a strongly $M$-Jacobi field along $\gamma_v$ implies a Jacobi 
field $Y$ along $\gamma_v$ satifying $Y(0)\in T_xM$ and $Y'(0)=-A_vY(0)$, 
where $Y'(0)=\widetilde{\nabla}_vY$ and $Y'(0)_T$ is the tangential (to $M$) 
component of $Y'(0)$.  
We call ${\rm Span}\{Y_0\,\vert\,Y\,:\,{\rm a}\,\,M{\rm -Jacobi}\,\,
{\rm field}\,\,{\rm along}\,\,\gamma_v\,\,{\rm s.t.}\,\,\lim_{t\to\infty}
\frac{\vert\vert Y_t\vert\vert}{t}=0\}$ the {\it nullity space} of the focal 
point $\gamma_v(\infty)$.  
Also, if there exists a $M$-Jacobi field $Y$ along $\gamma_v$ satisfying 
$\lim\limits_{t\to\infty}\dfrac{\vert\vert Y_t\vert\vert}{t}=0$ and 
${\rm Sec}(v,Y(0))<0$, then we call $\gamma_v(\infty)$ a {\it focal point of 
non-Euclidean type on} $N(\infty)$ {\it of} $M$ {\it along} $\gamma_v$, where 
${\rm Sec}(v,Y(0))$ is the sectional curvature for the $2$-plane spanned by 
$v$ and $Y(0)$.  
If $\exp^{\perp}(T^{\perp}_xM)$ is totally geodesic for each $x\in M$, $M$ is 
called a {\it submanifold with section}.  This notion has been recently 
defined in [HLO].  For a submanifold with section in a symmetric 
space of non-compact type, we have the following fact.  

\vspace{0.5truecm}

\centerline{
\unitlength 0.1in
\begin{picture}( 26.1000, 24.0000)( 22.0000,-34.0000)
%
\special{pn 8}%
\special{ar 3400 2200 1200 1200  0.0000000 0.0100000}%
\special{ar 3400 2200 1200 1200  0.0400000 0.0500000}%
\special{ar 3400 2200 1200 1200  0.0800000 0.0900000}%
\special{ar 3400 2200 1200 1200  0.1200000 0.1300000}%
\special{ar 3400 2200 1200 1200  0.1600000 0.1700000}%
\special{ar 3400 2200 1200 1200  0.2000000 0.2100000}%
\special{ar 3400 2200 1200 1200  0.2400000 0.2500000}%
\special{ar 3400 2200 1200 1200  0.2800000 0.2900000}%
\special{ar 3400 2200 1200 1200  0.3200000 0.3300000}%
\special{ar 3400 2200 1200 1200  0.3600000 0.3700000}%
\special{ar 3400 2200 1200 1200  0.4000000 0.4100000}%
\special{ar 3400 2200 1200 1200  0.4400000 0.4500000}%
\special{ar 3400 2200 1200 1200  0.4800000 0.4900000}%
\special{ar 3400 2200 1200 1200  0.5200000 0.5300000}%
\special{ar 3400 2200 1200 1200  0.5600000 0.5700000}%
\special{ar 3400 2200 1200 1200  0.6000000 0.6100000}%
\special{ar 3400 2200 1200 1200  0.6400000 0.6500000}%
\special{ar 3400 2200 1200 1200  0.6800000 0.6900000}%
\special{ar 3400 2200 1200 1200  0.7200000 0.7300000}%
\special{ar 3400 2200 1200 1200  0.7600000 0.7700000}%
\special{ar 3400 2200 1200 1200  0.8000000 0.8100000}%
\special{ar 3400 2200 1200 1200  0.8400000 0.8500000}%
\special{ar 3400 2200 1200 1200  0.8800000 0.8900000}%
\special{ar 3400 2200 1200 1200  0.9200000 0.9300000}%
\special{ar 3400 2200 1200 1200  0.9600000 0.9700000}%
\special{ar 3400 2200 1200 1200  1.0000000 1.0100000}%
\special{ar 3400 2200 1200 1200  1.0400000 1.0500000}%
\special{ar 3400 2200 1200 1200  1.0800000 1.0900000}%
\special{ar 3400 2200 1200 1200  1.1200000 1.1300000}%
\special{ar 3400 2200 1200 1200  1.1600000 1.1700000}%
\special{ar 3400 2200 1200 1200  1.2000000 1.2100000}%
\special{ar 3400 2200 1200 1200  1.2400000 1.2500000}%
\special{ar 3400 2200 1200 1200  1.2800000 1.2900000}%
\special{ar 3400 2200 1200 1200  1.3200000 1.3300000}%
\special{ar 3400 2200 1200 1200  1.3600000 1.3700000}%
\special{ar 3400 2200 1200 1200  1.4000000 1.4100000}%
\special{ar 3400 2200 1200 1200  1.4400000 1.4500000}%
\special{ar 3400 2200 1200 1200  1.4800000 1.4900000}%
\special{ar 3400 2200 1200 1200  1.5200000 1.5300000}%
\special{ar 3400 2200 1200 1200  1.5600000 1.5700000}%
\special{ar 3400 2200 1200 1200  1.6000000 1.6100000}%
\special{ar 3400 2200 1200 1200  1.6400000 1.6500000}%
\special{ar 3400 2200 1200 1200  1.6800000 1.6900000}%
\special{ar 3400 2200 1200 1200  1.7200000 1.7300000}%
\special{ar 3400 2200 1200 1200  1.7600000 1.7700000}%
\special{ar 3400 2200 1200 1200  1.8000000 1.8100000}%
\special{ar 3400 2200 1200 1200  1.8400000 1.8500000}%
\special{ar 3400 2200 1200 1200  1.8800000 1.8900000}%
\special{ar 3400 2200 1200 1200  1.9200000 1.9300000}%
\special{ar 3400 2200 1200 1200  1.9600000 1.9700000}%
\special{ar 3400 2200 1200 1200  2.0000000 2.0100000}%
\special{ar 3400 2200 1200 1200  2.0400000 2.0500000}%
\special{ar 3400 2200 1200 1200  2.0800000 2.0900000}%
\special{ar 3400 2200 1200 1200  2.1200000 2.1300000}%
\special{ar 3400 2200 1200 1200  2.1600000 2.1700000}%
\special{ar 3400 2200 1200 1200  2.2000000 2.2100000}%
\special{ar 3400 2200 1200 1200  2.2400000 2.2500000}%
\special{ar 3400 2200 1200 1200  2.2800000 2.2900000}%
\special{ar 3400 2200 1200 1200  2.3200000 2.3300000}%
\special{ar 3400 2200 1200 1200  2.3600000 2.3700000}%
\special{ar 3400 2200 1200 1200  2.4000000 2.4100000}%
\special{ar 3400 2200 1200 1200  2.4400000 2.4500000}%
\special{ar 3400 2200 1200 1200  2.4800000 2.4900000}%
\special{ar 3400 2200 1200 1200  2.5200000 2.5300000}%
\special{ar 3400 2200 1200 1200  2.5600000 2.5700000}%
\special{ar 3400 2200 1200 1200  2.6000000 2.6100000}%
\special{ar 3400 2200 1200 1200  2.6400000 2.6500000}%
\special{ar 3400 2200 1200 1200  2.6800000 2.6900000}%
\special{ar 3400 2200 1200 1200  2.7200000 2.7300000}%
\special{ar 3400 2200 1200 1200  2.7600000 2.7700000}%
\special{ar 3400 2200 1200 1200  2.8000000 2.8100000}%
\special{ar 3400 2200 1200 1200  2.8400000 2.8500000}%
\special{ar 3400 2200 1200 1200  2.8800000 2.8900000}%
\special{ar 3400 2200 1200 1200  2.9200000 2.9300000}%
\special{ar 3400 2200 1200 1200  2.9600000 2.9700000}%
\special{ar 3400 2200 1200 1200  3.0000000 3.0100000}%
\special{ar 3400 2200 1200 1200  3.0400000 3.0500000}%
\special{ar 3400 2200 1200 1200  3.0800000 3.0900000}%
\special{ar 3400 2200 1200 1200  3.1200000 3.1300000}%
\special{ar 3400 2200 1200 1200  3.1600000 3.1700000}%
\special{ar 3400 2200 1200 1200  3.2000000 3.2100000}%
\special{ar 3400 2200 1200 1200  3.2400000 3.2500000}%
\special{ar 3400 2200 1200 1200  3.2800000 3.2900000}%
\special{ar 3400 2200 1200 1200  3.3200000 3.3300000}%
\special{ar 3400 2200 1200 1200  3.3600000 3.3700000}%
\special{ar 3400 2200 1200 1200  3.4000000 3.4100000}%
\special{ar 3400 2200 1200 1200  3.4400000 3.4500000}%
\special{ar 3400 2200 1200 1200  3.4800000 3.4900000}%
\special{ar 3400 2200 1200 1200  3.5200000 3.5300000}%
\special{ar 3400 2200 1200 1200  3.5600000 3.5700000}%
\special{ar 3400 2200 1200 1200  3.6000000 3.6100000}%
\special{ar 3400 2200 1200 1200  3.6400000 3.6500000}%
\special{ar 3400 2200 1200 1200  3.6800000 3.6900000}%
\special{ar 3400 2200 1200 1200  3.7200000 3.7300000}%
\special{ar 3400 2200 1200 1200  3.7600000 3.7700000}%
\special{ar 3400 2200 1200 1200  3.8000000 3.8100000}%
\special{ar 3400 2200 1200 1200  3.8400000 3.8500000}%
\special{ar 3400 2200 1200 1200  3.8800000 3.8900000}%
\special{ar 3400 2200 1200 1200  3.9200000 3.9300000}%
\special{ar 3400 2200 1200 1200  3.9600000 3.9700000}%
\special{ar 3400 2200 1200 1200  4.0000000 4.0100000}%
\special{ar 3400 2200 1200 1200  4.0400000 4.0500000}%
\special{ar 3400 2200 1200 1200  4.0800000 4.0900000}%
\special{ar 3400 2200 1200 1200  4.1200000 4.1300000}%
\special{ar 3400 2200 1200 1200  4.1600000 4.1700000}%
\special{ar 3400 2200 1200 1200  4.2000000 4.2100000}%
\special{ar 3400 2200 1200 1200  4.2400000 4.2500000}%
\special{ar 3400 2200 1200 1200  4.2800000 4.2900000}%
\special{ar 3400 2200 1200 1200  4.3200000 4.3300000}%
\special{ar 3400 2200 1200 1200  4.3600000 4.3700000}%
\special{ar 3400 2200 1200 1200  4.4000000 4.4100000}%
\special{ar 3400 2200 1200 1200  4.4400000 4.4500000}%
\special{ar 3400 2200 1200 1200  4.4800000 4.4900000}%
\special{ar 3400 2200 1200 1200  4.5200000 4.5300000}%
\special{ar 3400 2200 1200 1200  4.5600000 4.5700000}%
\special{ar 3400 2200 1200 1200  4.6000000 4.6100000}%
\special{ar 3400 2200 1200 1200  4.6400000 4.6500000}%
\special{ar 3400 2200 1200 1200  4.6800000 4.6900000}%
\special{ar 3400 2200 1200 1200  4.7200000 4.7300000}%
\special{ar 3400 2200 1200 1200  4.7600000 4.7700000}%
\special{ar 3400 2200 1200 1200  4.8000000 4.8100000}%
\special{ar 3400 2200 1200 1200  4.8400000 4.8500000}%
\special{ar 3400 2200 1200 1200  4.8800000 4.8900000}%
\special{ar 3400 2200 1200 1200  4.9200000 4.9300000}%
\special{ar 3400 2200 1200 1200  4.9600000 4.9700000}%
\special{ar 3400 2200 1200 1200  5.0000000 5.0100000}%
\special{ar 3400 2200 1200 1200  5.0400000 5.0500000}%
\special{ar 3400 2200 1200 1200  5.0800000 5.0900000}%
\special{ar 3400 2200 1200 1200  5.1200000 5.1300000}%
\special{ar 3400 2200 1200 1200  5.1600000 5.1700000}%
\special{ar 3400 2200 1200 1200  5.2000000 5.2100000}%
\special{ar 3400 2200 1200 1200  5.2400000 5.2500000}%
\special{ar 3400 2200 1200 1200  5.2800000 5.2900000}%
\special{ar 3400 2200 1200 1200  5.3200000 5.3300000}%
\special{ar 3400 2200 1200 1200  5.3600000 5.3700000}%
\special{ar 3400 2200 1200 1200  5.4000000 5.4100000}%
\special{ar 3400 2200 1200 1200  5.4400000 5.4500000}%
\special{ar 3400 2200 1200 1200  5.4800000 5.4900000}%
\special{ar 3400 2200 1200 1200  5.5200000 5.5300000}%
\special{ar 3400 2200 1200 1200  5.5600000 5.5700000}%
\special{ar 3400 2200 1200 1200  5.6000000 5.6100000}%
\special{ar 3400 2200 1200 1200  5.6400000 5.6500000}%
\special{ar 3400 2200 1200 1200  5.6800000 5.6900000}%
\special{ar 3400 2200 1200 1200  5.7200000 5.7300000}%
\special{ar 3400 2200 1200 1200  5.7600000 5.7700000}%
\special{ar 3400 2200 1200 1200  5.8000000 5.8100000}%
\special{ar 3400 2200 1200 1200  5.8400000 5.8500000}%
\special{ar 3400 2200 1200 1200  5.8800000 5.8900000}%
\special{ar 3400 2200 1200 1200  5.9200000 5.9300000}%
\special{ar 3400 2200 1200 1200  5.9600000 5.9700000}%
\special{ar 3400 2200 1200 1200  6.0000000 6.0100000}%
\special{ar 3400 2200 1200 1200  6.0400000 6.0500000}%
\special{ar 3400 2200 1200 1200  6.0800000 6.0900000}%
\special{ar 3400 2200 1200 1200  6.1200000 6.1300000}%
\special{ar 3400 2200 1200 1200  6.1600000 6.1700000}%
\special{ar 3400 2200 1200 1200  6.2000000 6.2100000}%
\special{ar 3400 2200 1200 1200  6.2400000 6.2500000}%
\special{ar 3400 2200 1200 1200  6.2800000 6.2832853}%
%
\special{pn 8}%
\special{ar 3200 2200 570 610  2.0510899 4.2258370}%
%
\special{pn 8}%
\special{pa 2640 2200}%
\special{pa 4600 2200}%
\special{fp}%
\special{pa 4600 2200}%
\special{pa 4600 2200}%
\special{fp}%
%
\special{pn 8}%
\special{ar 5200 1600 2450 610  1.8178279 1.8256710}%
\special{ar 5200 1600 2450 610  1.8492005 1.8570436}%
\special{ar 5200 1600 2450 610  1.8805730 1.8884161}%
\special{ar 5200 1600 2450 610  1.9119456 1.9197887}%
\special{ar 5200 1600 2450 610  1.9433181 1.9511612}%
\special{ar 5200 1600 2450 610  1.9746907 1.9825338}%
\special{ar 5200 1600 2450 610  2.0060632 2.0139063}%
\special{ar 5200 1600 2450 610  2.0374357 2.0452789}%
\special{ar 5200 1600 2450 610  2.0688083 2.0766514}%
\special{ar 5200 1600 2450 610  2.1001808 2.1080240}%
\special{ar 5200 1600 2450 610  2.1315534 2.1393965}%
\special{ar 5200 1600 2450 610  2.1629259 2.1707691}%
\special{ar 5200 1600 2450 610  2.1942985 2.2021416}%
\special{ar 5200 1600 2450 610  2.2256710 2.2335142}%
\special{ar 5200 1600 2450 610  2.2570436 2.2648867}%
\special{ar 5200 1600 2450 610  2.2884161 2.2962593}%
\special{ar 5200 1600 2450 610  2.3197887 2.3276318}%
\special{ar 5200 1600 2450 610  2.3511612 2.3590044}%
\special{ar 5200 1600 2450 610  2.3825338 2.3903769}%
\special{ar 5200 1600 2450 610  2.4139063 2.4217495}%
\special{ar 5200 1600 2450 610  2.4452789 2.4531220}%
\special{ar 5200 1600 2450 610  2.4766514 2.4844946}%
\special{ar 5200 1600 2450 610  2.5080240 2.5158671}%
\special{ar 5200 1600 2450 610  2.5393965 2.5472397}%
\special{ar 5200 1600 2450 610  2.5707691 2.5786122}%
\special{ar 5200 1600 2450 610  2.6021416 2.6099848}%
\special{ar 5200 1600 2450 610  2.6335142 2.6413573}%
\special{ar 5200 1600 2450 610  2.6648867 2.6727299}%
\special{ar 5200 1600 2450 610  2.6962593 2.7041024}%
\special{ar 5200 1600 2450 610  2.7276318 2.7354750}%
\special{ar 5200 1600 2450 610  2.7590044 2.7668475}%
\special{ar 5200 1600 2450 610  2.7903769 2.7982201}%
\special{ar 5200 1600 2450 610  2.8217495 2.8295926}%
\special{ar 5200 1600 2450 610  2.8531220 2.8609652}%
%
\special{pn 8}%
\special{ar 5200 2800 2490 620  3.4194592 3.4271763}%
\special{ar 5200 2800 2490 620  3.4503274 3.4580444}%
\special{ar 5200 2800 2490 620  3.4811956 3.4889126}%
\special{ar 5200 2800 2490 620  3.5120637 3.5197808}%
\special{ar 5200 2800 2490 620  3.5429319 3.5506489}%
\special{ar 5200 2800 2490 620  3.5738001 3.5815171}%
\special{ar 5200 2800 2490 620  3.6046682 3.6123853}%
\special{ar 5200 2800 2490 620  3.6355364 3.6432534}%
\special{ar 5200 2800 2490 620  3.6664046 3.6741216}%
\special{ar 5200 2800 2490 620  3.6972727 3.7049898}%
\special{ar 5200 2800 2490 620  3.7281409 3.7358579}%
\special{ar 5200 2800 2490 620  3.7590091 3.7667261}%
\special{ar 5200 2800 2490 620  3.7898772 3.7975943}%
\special{ar 5200 2800 2490 620  3.8207454 3.8284624}%
\special{ar 5200 2800 2490 620  3.8516136 3.8593306}%
\special{ar 5200 2800 2490 620  3.8824817 3.8901988}%
\special{ar 5200 2800 2490 620  3.9133499 3.9210669}%
\special{ar 5200 2800 2490 620  3.9442181 3.9519351}%
\special{ar 5200 2800 2490 620  3.9750862 3.9828033}%
\special{ar 5200 2800 2490 620  4.0059544 4.0136714}%
\special{ar 5200 2800 2490 620  4.0368226 4.0445396}%
\special{ar 5200 2800 2490 620  4.0676907 4.0754078}%
\special{ar 5200 2800 2490 620  4.0985589 4.1062759}%
\special{ar 5200 2800 2490 620  4.1294271 4.1371441}%
\special{ar 5200 2800 2490 620  4.1602952 4.1680123}%
\special{ar 5200 2800 2490 620  4.1911634 4.1988804}%
\special{ar 5200 2800 2490 620  4.2220316 4.2297486}%
\special{ar 5200 2800 2490 620  4.2528997 4.2606168}%
\special{ar 5200 2800 2490 620  4.2837679 4.2914849}%
\special{ar 5200 2800 2490 620  4.3146361 4.3223531}%
\special{ar 5200 2800 2490 620  4.3455042 4.3532213}%
\special{ar 5200 2800 2490 620  4.3763724 4.3840894}%
\special{ar 5200 2800 2490 620  4.4072406 4.4149576}%
\special{ar 5200 2800 2490 620  4.4381087 4.4458258}%
%
\special{pn 8}%
\special{pa 2640 2200}%
\special{pa 2640 1710}%
\special{fp}%
\special{sh 1}%
\special{pa 2640 1710}%
\special{pa 2620 1778}%
\special{pa 2640 1764}%
\special{pa 2660 1778}%
\special{pa 2640 1710}%
\special{fp}%
%
\special{pn 8}%
\special{pa 2960 2200}%
\special{pa 2960 1850}%
\special{fp}%
\special{sh 1}%
\special{pa 2960 1850}%
\special{pa 2940 1918}%
\special{pa 2960 1904}%
\special{pa 2980 1918}%
\special{pa 2960 1850}%
\special{fp}%
%
\special{pn 8}%
\special{pa 3290 2200}%
\special{pa 3290 1980}%
\special{fp}%
\special{sh 1}%
\special{pa 3290 1980}%
\special{pa 3270 2048}%
\special{pa 3290 2034}%
\special{pa 3310 2048}%
\special{pa 3290 1980}%
\special{fp}%
%
\special{pn 8}%
\special{pa 3630 2200}%
\special{pa 3630 2060}%
\special{fp}%
\special{sh 1}%
\special{pa 3630 2060}%
\special{pa 3610 2128}%
\special{pa 3630 2114}%
\special{pa 3650 2128}%
\special{pa 3630 2060}%
\special{fp}%
%
\special{pn 8}%
\special{pa 3930 2200}%
\special{pa 3930 2120}%
\special{fp}%
\special{sh 1}%
\special{pa 3930 2120}%
\special{pa 3910 2188}%
\special{pa 3930 2174}%
\special{pa 3950 2188}%
\special{pa 3930 2120}%
\special{fp}%
%
\special{pn 8}%
\special{pa 3800 2590}%
\special{pa 3450 2200}%
\special{dt 0.045}%
\special{sh 1}%
\special{pa 3450 2200}%
\special{pa 3480 2264}%
\special{pa 3486 2240}%
\special{pa 3510 2236}%
\special{pa 3450 2200}%
\special{fp}%
%
\special{pn 13}%
\special{sh 1}%
\special{ar 4600 2200 10 10 0  6.28318530717959E+0000}%
\special{sh 1}%
\special{ar 4600 2200 10 10 0  6.28318530717959E+0000}%
%
\special{pn 8}%
\special{pa 4760 1330}%
\special{pa 4390 1530}%
\special{dt 0.045}%
\special{sh 1}%
\special{pa 4390 1530}%
\special{pa 4458 1516}%
\special{pa 4438 1506}%
\special{pa 4440 1482}%
\special{pa 4390 1530}%
\special{fp}%
%
\special{pn 8}%
\special{pa 4810 2050}%
\special{pa 4610 2200}%
\special{dt 0.045}%
\special{sh 1}%
\special{pa 4610 2200}%
\special{pa 4676 2176}%
\special{pa 4654 2168}%
\special{pa 4652 2144}%
\special{pa 4610 2200}%
\special{fp}%
\put(31.9000,-27.2000){\makebox(0,0)[rt]{$M$}}%
\put(39.5000,-26.3000){\makebox(0,0)[rt]{$\gamma_v$}}%
\put(35.0000,-17.5000){\makebox(0,0)[rt]{$Y$}}%
\put(51.2000,-11.5000){\makebox(0,0)[rt]{$N(\infty)$}}%
\put(50.8000,-17.9000){\makebox(0,0)[rt]{$\gamma_v(\infty)$}}%
\end{picture}%
}

\vspace{0.3truecm}

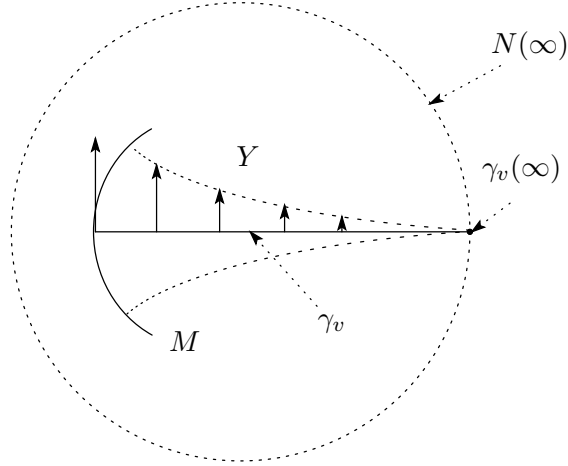
\centerline{{\bf Fig. 1.}}

\vspace{0.5truecm}

\noindent
{\bf Proposition 3.1.} {\sl Let $M$ be a submanifold with section in 
a symmetric space $N:=G/K$ of non-compact type and $v$ be a normal vector of 
$M$ at $x$.  Then the following conditions {\rm (i)} and {\rm (ii)} are 
equivalent:

{\rm (i)} $\gamma_v(\infty)$ is a focal point on $N(\infty)$ of $M$ along 
$\gamma_v$,

{\rm (ii)} $\gamma_v(\infty)$ is a strongly focal point on $N(\infty)$ of $M$ 
along $\gamma_v$.  

\noindent
Furthermore, if $M$ is homogeneous (hence it is a principal orbit of a polar 
action $H$ on $N$), then these conditions are equivalent to 
the following conditions:

{\rm (iii)} there exists a normal geodesic variation 
$\delta:[0,\infty)\times(-\varepsilon,\varepsilon)\to N$ such that 
$\delta(\cdot,0)=\gamma_v(\cdot)$, the variational vector field 
$\frac{\partial\delta}{\partial s}\vert_{s=0}$ is a strongly $M$-Jacobi 
field and that 
$\delta(\cdot,s)(\infty)=\gamma_v(\infty)$ for any $s\in(-\varepsilon,
\varepsilon)$.

{\rm (iv)} the action on $N(\infty)$ induced from the $H$-action posseses 
a non-trivial subaction having $\gamma_v(\infty)$ as a fixed point.}

\vspace{0.5truecm}

\noindent
{\it Proof.} First we shall show ${\rm (i)}\Rightarrow{\rm (ii)}$.  
Assume that $\gamma_v(\infty)$ is a focal point on $N(\infty)$ along 
$\gamma_v$.  Hence there exists an $M$-Jacobi field $Y$ along $\gamma_v$ 
such that $\lim\limits_{t\to\infty}\frac{\vert\vert Y(t)\vert\vert}{t}=0$.  
The Jacobi field $Y$ is described as 
$$Y(t)=P_{\gamma_v\vert_{[0,t]}}\left(D^{co}_{tv}(Y(0))+D^{si}_{tv}
(-A_{tv}Y(0)+Y'(0)_{\perp})\right),$$
where $P_{\gamma_v\vert_{[0,t]}}$ is the parallel translation along 
$\gamma_v\vert_{[0,t]}$, $D^{co}_{tv}$ and $D^{si}_{tv}$ are the above 
operators, $A$ (resp. $\nabla^{\perp}$) is the shape tensor (resp. the normal 
connection) of $M$ and $Y'(0)_{\perp}$ is the normal component of 
$Y'(0)(=\widetilde{\nabla}_vY)$.  Since $M$ has section, we have 
$D^{co}_{tv}(Y(0)),\,D^{si}_{tv}(A_{tv}Y(0))\in T_xM$.  Hence we have 
$\vert\vert Y(t)\vert\vert\geq\vert\vert(D^{co}_{tv}-D^{si}_{tv}\circ A_{tv})
(Y(0))\vert\vert$.  The strongly $M$-Jacobi field $Y^S$ along $\gamma_v$ with 
$Y^S(0)=Y(0)$ is described as 
$$Y^S(t)=P_{\gamma_v\vert_{[0,t]}}\left(
(D^{co}_{tv}-D^{si}_{tv}\circ A_{tv})(Y(0))\right). \leqno{(3.1)}$$
Hence we have $\vert\vert Y(t)\vert\vert\geq\vert\vert Y^S(t)\vert\vert$ and 
hence $\lim\limits_{t\to 0}\frac{\vert\vert Y^S(t)\vert\vert}{t}=0$.  Thus 
$\gamma_v(\infty)$ is a strongly focal point on $N(\infty)$ along 
$\gamma_v$.  Thus we have ${\rm (i)}\Rightarrow{\rm (ii)}$.  The converse 
${\rm (ii)}\Rightarrow{\rm (i)}$ is trivial.  

Next we shall show that ${\rm (ii)}\Rightarrow{\rm (iii)}$ holds if $M$ is 
homogeneous.  Assume that $\gamma_v(\infty)$ is a strongly focal point 
on $N(\infty)$ along $\gamma_v$.  Hence there exists a strongly $M$-Jacobi 
field $Y^S$ along $\gamma_v$ with $\lim\limits_{t\to 0}\frac{\vert\vert Y^S(t)
\vert\vert}{t}=0$.  Since $Y^S$ is described as in $(3.1)$, we have 
$\vert\vert Y^S(t)\vert\vert=\vert\vert(D^{co}_{tv}-D^{si}_{tv}\circ A_{tv})
(Y^S(0))\vert\vert$.  
Since $M$ is a homogeneous submanifold with section, it is catched as a 
principal orbit of some complex polar action $H\times G/K\to G/K\,\,
(H\subset G)$.  See [Koi2] about the definition of a complex polar action.  
Let $\{\exp\,sX\,\vert\,s\in{\bf R}\}$ be a one-parameter subgroup of $H$ such 
that $\frac{d(\exp\,sX)(x)}{ds}\vert_{s=0}=Y(0)$.  Set $\alpha(s):=
(\exp\,sX)(x)$.  Let $\widetilde v$ be the parallel normal vector field along 
$\alpha$ with $\widetilde v_0=v$.  Define $\delta:[0,\infty)\times
(-\varepsilon,\varepsilon)\to N$ by $\delta(t,s):=\exp^{\perp}
(t\widetilde v_s)$.  We have $\frac{\partial\delta}{\partial s}\vert_{s=0}
=Y^S$.  Set $Y^S_{s_0}:=\frac{\partial\delta}{\partial s}
\vert_{[0,\infty)\times\{s_0\}}$ for each $s_0\in(-\varepsilon,\varepsilon)$.  
Since $Y^S_{s_0}$ is a strongly $M$-Jacobi field along $\delta(\cdot,s_0)$, 
it is described as in $(3.1)$.  Hence we have $\vert\vert Y^S_{s_0}(t)
\vert\vert=\vert\vert(D^{co}_{t\widetilde v_s}-D^{si}_{t\widetilde v_s}\circ 
A_{t\widetilde v_s})(Y^S_{s_0}(0))\vert\vert$.  Since $M$ is a principal orbit 
of the $H$-action, we have $\widetilde v_s=(\exp\,sX)_{\ast}(v)$ 
($s\in(-\varepsilon,\varepsilon)$).  From this fact, we have 
$\vert\vert(D^{co}_{t\widetilde v_{s_0}}-D^{si}_{t\widetilde v_{s_0}}\circ 
A_{t\widetilde v_{s_0}})(Y^S_{s_0}(0))\vert\vert=
\vert\vert(D^{co}_{tv}-D^{si}_{tv}\circ A_{tv})(Y^S(0))\vert\vert$ 
($(t,s_0)\in[0,\infty)\times(-\varepsilon,\varepsilon)$).  Therefore, we have 
$\vert\vert Y^S_{s_0}(t)\vert\vert=\vert\vert Y^S(t)\vert\vert$ 
($(t,s_0)\in[0,\infty)\times(-\varepsilon,\varepsilon)$).  Hence we have 
$$\begin{array}{l}
\displaystyle{\lim_{t\to\infty}\frac{d(\delta(t,s_0),\gamma_v(t))}{t}\leq
\lim_{t\to\infty}\frac 1t\int_0^{s_0}\vert\vert Y^S_s(t)\vert\vert ds}\\
\hspace{3.5truecm}\displaystyle{=\lim_{t\to\infty}\frac{s_0\vert\vert 
Y^S(t)\vert\vert}{t}=0,}
\end{array}$$
that is, $\delta(\cdot,s_0)(\infty)=\gamma_v(\infty)$.  
Thus ${\rm (ii)}\Rightarrow{\rm (iii)}$ is shown.  The converse is trivial.  
Also, ${\rm (iii)}\Leftrightarrow{\rm (iv)}$ is trivial.  
This completes the proof.  \hspace{6truecm}q.e.d.

\vspace{0.5truecm}

\noindent
{\it Remark 3.1.} Let $\gamma$ be a normal geodesic of a princial orbit $M$ 
of a polar action $H$ on $N=G/K$.  If $\gamma(\infty)$ is a fixed point of 
the action on $N(\infty)$ induced from the $H$-action, then $\gamma(\infty)$ 
is a focal point of $M$ along $\gamma$ having $T_{\gamma(0)}M$ as the nullity 
space.  

\vspace{0.5truecm}

Now we shall illustrate that the second-half of the statement in 
Proposition 3.1 does not hold without the assumption of the homogeneity of 
$M$.  Let $S$ be a horosphere in a symmetric space $N=G/K$ of non-compact 
type and $M$ be a non-homogeneous hypersurface in $N$ through $x\in S$ such 
that $j^2_x(\iota_M)=j^2_x(\iota_S)$ but $j^3_x(\iota_M)\not=j^3_x(\iota_S)$ 
and that $M$ positions outside or inside $S$ (see Fig. 2), where $\iota_M$ 
(resp. $\iota_S$) is the inclusion map of $M$ (resp. $S$) into $N$ and 
$j^2_x(\cdot)$ is the $2$-jet of $\cdot$ at $x$.  Since $M$ is a hypersurface, 
it has sections.  Let $v$ be the inward unit normal vector of $S$ at $x$.  
Then $\gamma_v(\infty)$ is a focal point on $N(\infty)$ of $M$ along 
$\gamma_v$ but there does not exist a normal (to $M$) geodesic variation 
$\delta:[0,\infty)\times(-\varepsilon,\varepsilon)\to N$ of $\gamma_v$ 
such that $\delta(\cdot,0)=\gamma_v$, the variational vector field 
$\frac{\partial\delta}{\partial s}\vert_{s=0}$ is a strongly $M$-Jacobi field 
and that $\delta(\cdot,s_0)(\infty)=\gamma_v(\infty)$ for each 
$s_0\in(-\varepsilon,\varepsilon)$.  
Thus the second-half of the statement in Proposition 3.1 does not hold 
without the assumption of the homogeneity of $M$.  
Since $S$ is complex equifocal, $j^2_x(\iota_M)=j^2_x(\iota_S)$ and 
$j^3_x(\iota_M)\not=j^3_x(\iota_S)$, we see that $M$ is not 
complex equifocal.  In more general, it is conjectured that the second-half 
of the statement in Proposition 3.1 holds if $M$ is complex equifocal.  

\vspace{0.5truecm}

\centerline{
\unitlength 0.1in
\begin{picture}( 27.9000, 24.0000)( 20.2000,-34.0000)
%
\special{pn 8}%
\special{pa 2640 2200}%
\special{pa 4600 2200}%
\special{fp}%
\special{pa 4600 2200}%
\special{pa 4600 2200}%
\special{fp}%
%
\special{pn 8}%
\special{ar 5200 1600 2450 610  1.8178279 1.8256710}%
\special{ar 5200 1600 2450 610  1.8492005 1.8570436}%
\special{ar 5200 1600 2450 610  1.8805730 1.8884161}%
\special{ar 5200 1600 2450 610  1.9119456 1.9197887}%
\special{ar 5200 1600 2450 610  1.9433181 1.9511612}%
\special{ar 5200 1600 2450 610  1.9746907 1.9825338}%
\special{ar 5200 1600 2450 610  2.0060632 2.0139063}%
\special{ar 5200 1600 2450 610  2.0374357 2.0452789}%
\special{ar 5200 1600 2450 610  2.0688083 2.0766514}%
\special{ar 5200 1600 2450 610  2.1001808 2.1080240}%
\special{ar 5200 1600 2450 610  2.1315534 2.1393965}%
\special{ar 5200 1600 2450 610  2.1629259 2.1707691}%
\special{ar 5200 1600 2450 610  2.1942985 2.2021416}%
\special{ar 5200 1600 2450 610  2.2256710 2.2335142}%
\special{ar 5200 1600 2450 610  2.2570436 2.2648867}%
\special{ar 5200 1600 2450 610  2.2884161 2.2962593}%
\special{ar 5200 1600 2450 610  2.3197887 2.3276318}%
\special{ar 5200 1600 2450 610  2.3511612 2.3590044}%
\special{ar 5200 1600 2450 610  2.3825338 2.3903769}%
\special{ar 5200 1600 2450 610  2.4139063 2.4217495}%
\special{ar 5200 1600 2450 610  2.4452789 2.4531220}%
\special{ar 5200 1600 2450 610  2.4766514 2.4844946}%
\special{ar 5200 1600 2450 610  2.5080240 2.5158671}%
\special{ar 5200 1600 2450 610  2.5393965 2.5472397}%
\special{ar 5200 1600 2450 610  2.5707691 2.5786122}%
\special{ar 5200 1600 2450 610  2.6021416 2.6099848}%
\special{ar 5200 1600 2450 610  2.6335142 2.6413573}%
\special{ar 5200 1600 2450 610  2.6648867 2.6727299}%
\special{ar 5200 1600 2450 610  2.6962593 2.7041024}%
\special{ar 5200 1600 2450 610  2.7276318 2.7354750}%
\special{ar 5200 1600 2450 610  2.7590044 2.7668475}%
\special{ar 5200 1600 2450 610  2.7903769 2.7982201}%
\special{ar 5200 1600 2450 610  2.8217495 2.8295926}%
\special{ar 5200 1600 2450 610  2.8531220 2.8609652}%
%
\special{pn 8}%
\special{ar 5200 2800 2490 620  3.4194592 3.4271763}%
\special{ar 5200 2800 2490 620  3.4503274 3.4580444}%
\special{ar 5200 2800 2490 620  3.4811956 3.4889126}%
\special{ar 5200 2800 2490 620  3.5120637 3.5197808}%
\special{ar 5200 2800 2490 620  3.5429319 3.5506489}%
\special{ar 5200 2800 2490 620  3.5738001 3.5815171}%
\special{ar 5200 2800 2490 620  3.6046682 3.6123853}%
\special{ar 5200 2800 2490 620  3.6355364 3.6432534}%
\special{ar 5200 2800 2490 620  3.6664046 3.6741216}%
\special{ar 5200 2800 2490 620  3.6972727 3.7049898}%
\special{ar 5200 2800 2490 620  3.7281409 3.7358579}%
\special{ar 5200 2800 2490 620  3.7590091 3.7667261}%
\special{ar 5200 2800 2490 620  3.7898772 3.7975943}%
\special{ar 5200 2800 2490 620  3.8207454 3.8284624}%
\special{ar 5200 2800 2490 620  3.8516136 3.8593306}%
\special{ar 5200 2800 2490 620  3.8824817 3.8901988}%
\special{ar 5200 2800 2490 620  3.9133499 3.9210669}%
\special{ar 5200 2800 2490 620  3.9442181 3.9519351}%
\special{ar 5200 2800 2490 620  3.9750862 3.9828033}%
\special{ar 5200 2800 2490 620  4.0059544 4.0136714}%
\special{ar 5200 2800 2490 620  4.0368226 4.0445396}%
\special{ar 5200 2800 2490 620  4.0676907 4.0754078}%
\special{ar 5200 2800 2490 620  4.0985589 4.1062759}%
\special{ar 5200 2800 2490 620  4.1294271 4.1371441}%
\special{ar 5200 2800 2490 620  4.1602952 4.1680123}%
\special{ar 5200 2800 2490 620  4.1911634 4.1988804}%
\special{ar 5200 2800 2490 620  4.2220316 4.2297486}%
\special{ar 5200 2800 2490 620  4.2528997 4.2606168}%
\special{ar 5200 2800 2490 620  4.2837679 4.2914849}%
\special{ar 5200 2800 2490 620  4.3146361 4.3223531}%
\special{ar 5200 2800 2490 620  4.3455042 4.3532213}%
\special{ar 5200 2800 2490 620  4.3763724 4.3840894}%
\special{ar 5200 2800 2490 620  4.4072406 4.4149576}%
\special{ar 5200 2800 2490 620  4.4381087 4.4458258}%
%
\special{pn 8}%
\special{pa 2640 2200}%
\special{pa 2640 1710}%
\special{fp}%
\special{sh 1}%
\special{pa 2640 1710}%
\special{pa 2620 1778}%
\special{pa 2640 1764}%
\special{pa 2660 1778}%
\special{pa 2640 1710}%
\special{fp}%
%
\special{pn 8}%
\special{pa 2960 2200}%
\special{pa 2960 1850}%
\special{fp}%
\special{sh 1}%
\special{pa 2960 1850}%
\special{pa 2940 1918}%
\special{pa 2960 1904}%
\special{pa 2980 1918}%
\special{pa 2960 1850}%
\special{fp}%
%
\special{pn 8}%
\special{pa 3290 2200}%
\special{pa 3290 1980}%
\special{fp}%
\special{sh 1}%
\special{pa 3290 1980}%
\special{pa 3270 2048}%
\special{pa 3290 2034}%
\special{pa 3310 2048}%
\special{pa 3290 1980}%
\special{fp}%
%
\special{pn 8}%
\special{pa 3630 2200}%
\special{pa 3630 2060}%
\special{fp}%
\special{sh 1}%
\special{pa 3630 2060}%
\special{pa 3610 2128}%
\special{pa 3630 2114}%
\special{pa 3650 2128}%
\special{pa 3630 2060}%
\special{fp}%
%
\special{pn 8}%
\special{pa 3930 2200}%
\special{pa 3930 2120}%
\special{fp}%
\special{sh 1}%
\special{pa 3930 2120}%
\special{pa 3910 2188}%
\special{pa 3930 2174}%
\special{pa 3950 2188}%
\special{pa 3930 2120}%
\special{fp}%
%
\special{pn 13}%
\special{sh 1}%
\special{ar 4600 2200 10 10 0  6.28318530717959E+0000}%
\special{sh 1}%
\special{ar 4600 2200 10 10 0  6.28318530717959E+0000}%
%
\special{pn 8}%
\special{pa 4760 1330}%
\special{pa 4390 1530}%
\special{dt 0.045}%
\special{sh 1}%
\special{pa 4390 1530}%
\special{pa 4458 1516}%
\special{pa 4438 1506}%
\special{pa 4440 1482}%
\special{pa 4390 1530}%
\special{fp}%
%
\special{pn 8}%
\special{pa 4810 2050}%
\special{pa 4610 2200}%
\special{dt 0.045}%
\special{sh 1}%
\special{pa 4610 2200}%
\special{pa 4676 2176}%
\special{pa 4654 2168}%
\special{pa 4652 2144}%
\special{pa 4610 2200}%
\special{fp}%
\put(22.9000,-30.2000){\makebox(0,0)[rt]{$M$}}%
\put(42.7000,-28.6000){\makebox(0,0)[rt]{$\gamma_v$}}%
\put(51.2000,-11.5000){\makebox(0,0)[rt]{$N(\infty)$}}%
\put(50.8000,-17.9000){\makebox(0,0)[rt]{$\gamma_v(\infty)$}}%
%
\special{pn 8}%
\special{pa 4040 2810}%
\special{pa 3750 2200}%
\special{dt 0.045}%
\special{sh 1}%
\special{pa 3750 2200}%
\special{pa 3762 2270}%
\special{pa 3774 2248}%
\special{pa 3798 2252}%
\special{pa 3750 2200}%
\special{fp}%
\special{pa 3640 3020}%
\special{pa 3520 2650}%
\special{dt 0.045}%
\special{sh 1}%
\special{pa 3520 2650}%
\special{pa 3522 2720}%
\special{pa 3536 2702}%
\special{pa 3560 2708}%
\special{pa 3520 2650}%
\special{fp}%
\put(37.5000,-30.7000){\makebox(0,0)[rt]{$S$}}%
%
\special{pn 8}%
\special{ar 3400 2200 1200 1200  0.0000000 0.0100000}%
\special{ar 3400 2200 1200 1200  0.0400000 0.0500000}%
\special{ar 3400 2200 1200 1200  0.0800000 0.0900000}%
\special{ar 3400 2200 1200 1200  0.1200000 0.1300000}%
\special{ar 3400 2200 1200 1200  0.1600000 0.1700000}%
\special{ar 3400 2200 1200 1200  0.2000000 0.2100000}%
\special{ar 3400 2200 1200 1200  0.2400000 0.2500000}%
\special{ar 3400 2200 1200 1200  0.2800000 0.2900000}%
\special{ar 3400 2200 1200 1200  0.3200000 0.3300000}%
\special{ar 3400 2200 1200 1200  0.3600000 0.3700000}%
\special{ar 3400 2200 1200 1200  0.4000000 0.4100000}%
\special{ar 3400 2200 1200 1200  0.4400000 0.4500000}%
\special{ar 3400 2200 1200 1200  0.4800000 0.4900000}%
\special{ar 3400 2200 1200 1200  0.5200000 0.5300000}%
\special{ar 3400 2200 1200 1200  0.5600000 0.5700000}%
\special{ar 3400 2200 1200 1200  0.6000000 0.6100000}%
\special{ar 3400 2200 1200 1200  0.6400000 0.6500000}%
\special{ar 3400 2200 1200 1200  0.6800000 0.6900000}%
\special{ar 3400 2200 1200 1200  0.7200000 0.7300000}%
\special{ar 3400 2200 1200 1200  0.7600000 0.7700000}%
\special{ar 3400 2200 1200 1200  0.8000000 0.8100000}%
\special{ar 3400 2200 1200 1200  0.8400000 0.8500000}%
\special{ar 3400 2200 1200 1200  0.8800000 0.8900000}%
\special{ar 3400 2200 1200 1200  0.9200000 0.9300000}%
\special{ar 3400 2200 1200 1200  0.9600000 0.9700000}%
\special{ar 3400 2200 1200 1200  1.0000000 1.0100000}%
\special{ar 3400 2200 1200 1200  1.0400000 1.0500000}%
\special{ar 3400 2200 1200 1200  1.0800000 1.0900000}%
\special{ar 3400 2200 1200 1200  1.1200000 1.1300000}%
\special{ar 3400 2200 1200 1200  1.1600000 1.1700000}%
\special{ar 3400 2200 1200 1200  1.2000000 1.2100000}%
\special{ar 3400 2200 1200 1200  1.2400000 1.2500000}%
\special{ar 3400 2200 1200 1200  1.2800000 1.2900000}%
\special{ar 3400 2200 1200 1200  1.3200000 1.3300000}%
\special{ar 3400 2200 1200 1200  1.3600000 1.3700000}%
\special{ar 3400 2200 1200 1200  1.4000000 1.4100000}%
\special{ar 3400 2200 1200 1200  1.4400000 1.4500000}%
\special{ar 3400 2200 1200 1200  1.4800000 1.4900000}%
\special{ar 3400 2200 1200 1200  1.5200000 1.5300000}%
\special{ar 3400 2200 1200 1200  1.5600000 1.5700000}%
\special{ar 3400 2200 1200 1200  1.6000000 1.6100000}%
\special{ar 3400 2200 1200 1200  1.6400000 1.6500000}%
\special{ar 3400 2200 1200 1200  1.6800000 1.6900000}%
\special{ar 3400 2200 1200 1200  1.7200000 1.7300000}%
\special{ar 3400 2200 1200 1200  1.7600000 1.7700000}%
\special{ar 3400 2200 1200 1200  1.8000000 1.8100000}%
\special{ar 3400 2200 1200 1200  1.8400000 1.8500000}%
\special{ar 3400 2200 1200 1200  1.8800000 1.8900000}%
\special{ar 3400 2200 1200 1200  1.9200000 1.9300000}%
\special{ar 3400 2200 1200 1200  1.9600000 1.9700000}%
\special{ar 3400 2200 1200 1200  2.0000000 2.0100000}%
\special{ar 3400 2200 1200 1200  2.0400000 2.0500000}%
\special{ar 3400 2200 1200 1200  2.0800000 2.0900000}%
\special{ar 3400 2200 1200 1200  2.1200000 2.1300000}%
\special{ar 3400 2200 1200 1200  2.1600000 2.1700000}%
\special{ar 3400 2200 1200 1200  2.2000000 2.2100000}%
\special{ar 3400 2200 1200 1200  2.2400000 2.2500000}%
\special{ar 3400 2200 1200 1200  2.2800000 2.2900000}%
\special{ar 3400 2200 1200 1200  2.3200000 2.3300000}%
\special{ar 3400 2200 1200 1200  2.3600000 2.3700000}%
\special{ar 3400 2200 1200 1200  2.4000000 2.4100000}%
\special{ar 3400 2200 1200 1200  2.4400000 2.4500000}%
\special{ar 3400 2200 1200 1200  2.4800000 2.4900000}%
\special{ar 3400 2200 1200 1200  2.5200000 2.5300000}%
\special{ar 3400 2200 1200 1200  2.5600000 2.5700000}%
\special{ar 3400 2200 1200 1200  2.6000000 2.6100000}%
\special{ar 3400 2200 1200 1200  2.6400000 2.6500000}%
\special{ar 3400 2200 1200 1200  2.6800000 2.6900000}%
\special{ar 3400 2200 1200 1200  2.7200000 2.7300000}%
\special{ar 3400 2200 1200 1200  2.7600000 2.7700000}%
\special{ar 3400 2200 1200 1200  2.8000000 2.8100000}%
\special{ar 3400 2200 1200 1200  2.8400000 2.8500000}%
\special{ar 3400 2200 1200 1200  2.8800000 2.8900000}%
\special{ar 3400 2200 1200 1200  2.9200000 2.9300000}%
\special{ar 3400 2200 1200 1200  2.9600000 2.9700000}%
\special{ar 3400 2200 1200 1200  3.0000000 3.0100000}%
\special{ar 3400 2200 1200 1200  3.0400000 3.0500000}%
\special{ar 3400 2200 1200 1200  3.0800000 3.0900000}%
\special{ar 3400 2200 1200 1200  3.1200000 3.1300000}%
\special{ar 3400 2200 1200 1200  3.1600000 3.1700000}%
\special{ar 3400 2200 1200 1200  3.2000000 3.2100000}%
\special{ar 3400 2200 1200 1200  3.2400000 3.2500000}%
\special{ar 3400 2200 1200 1200  3.2800000 3.2900000}%
\special{ar 3400 2200 1200 1200  3.3200000 3.3300000}%
\special{ar 3400 2200 1200 1200  3.3600000 3.3700000}%
\special{ar 3400 2200 1200 1200  3.4000000 3.4100000}%
\special{ar 3400 2200 1200 1200  3.4400000 3.4500000}%
\special{ar 3400 2200 1200 1200  3.4800000 3.4900000}%
\special{ar 3400 2200 1200 1200  3.5200000 3.5300000}%
\special{ar 3400 2200 1200 1200  3.5600000 3.5700000}%
\special{ar 3400 2200 1200 1200  3.6000000 3.6100000}%
\special{ar 3400 2200 1200 1200  3.6400000 3.6500000}%
\special{ar 3400 2200 1200 1200  3.6800000 3.6900000}%
\special{ar 3400 2200 1200 1200  3.7200000 3.7300000}%
\special{ar 3400 2200 1200 1200  3.7600000 3.7700000}%
\special{ar 3400 2200 1200 1200  3.8000000 3.8100000}%
\special{ar 3400 2200 1200 1200  3.8400000 3.8500000}%
\special{ar 3400 2200 1200 1200  3.8800000 3.8900000}%
\special{ar 3400 2200 1200 1200  3.9200000 3.9300000}%
\special{ar 3400 2200 1200 1200  3.9600000 3.9700000}%
\special{ar 3400 2200 1200 1200  4.0000000 4.0100000}%
\special{ar 3400 2200 1200 1200  4.0400000 4.0500000}%
\special{ar 3400 2200 1200 1200  4.0800000 4.0900000}%
\special{ar 3400 2200 1200 1200  4.1200000 4.1300000}%
\special{ar 3400 2200 1200 1200  4.1600000 4.1700000}%
\special{ar 3400 2200 1200 1200  4.2000000 4.2100000}%
\special{ar 3400 2200 1200 1200  4.2400000 4.2500000}%
\special{ar 3400 2200 1200 1200  4.2800000 4.2900000}%
\special{ar 3400 2200 1200 1200  4.3200000 4.3300000}%
\special{ar 3400 2200 1200 1200  4.3600000 4.3700000}%
\special{ar 3400 2200 1200 1200  4.4000000 4.4100000}%
\special{ar 3400 2200 1200 1200  4.4400000 4.4500000}%
\special{ar 3400 2200 1200 1200  4.4800000 4.4900000}%
\special{ar 3400 2200 1200 1200  4.5200000 4.5300000}%
\special{ar 3400 2200 1200 1200  4.5600000 4.5700000}%
\special{ar 3400 2200 1200 1200  4.6000000 4.6100000}%
\special{ar 3400 2200 1200 1200  4.6400000 4.6500000}%
\special{ar 3400 2200 1200 1200  4.6800000 4.6900000}%
\special{ar 3400 2200 1200 1200  4.7200000 4.7300000}%
\special{ar 3400 2200 1200 1200  4.7600000 4.7700000}%
\special{ar 3400 2200 1200 1200  4.8000000 4.8100000}%
\special{ar 3400 2200 1200 1200  4.8400000 4.8500000}%
\special{ar 3400 2200 1200 1200  4.8800000 4.8900000}%
\special{ar 3400 2200 1200 1200  4.9200000 4.9300000}%
\special{ar 3400 2200 1200 1200  4.9600000 4.9700000}%
\special{ar 3400 2200 1200 1200  5.0000000 5.0100000}%
\special{ar 3400 2200 1200 1200  5.0400000 5.0500000}%
\special{ar 3400 2200 1200 1200  5.0800000 5.0900000}%
\special{ar 3400 2200 1200 1200  5.1200000 5.1300000}%
\special{ar 3400 2200 1200 1200  5.1600000 5.1700000}%
\special{ar 3400 2200 1200 1200  5.2000000 5.2100000}%
\special{ar 3400 2200 1200 1200  5.2400000 5.2500000}%
\special{ar 3400 2200 1200 1200  5.2800000 5.2900000}%
\special{ar 3400 2200 1200 1200  5.3200000 5.3300000}%
\special{ar 3400 2200 1200 1200  5.3600000 5.3700000}%
\special{ar 3400 2200 1200 1200  5.4000000 5.4100000}%
\special{ar 3400 2200 1200 1200  5.4400000 5.4500000}%
\special{ar 3400 2200 1200 1200  5.4800000 5.4900000}%
\special{ar 3400 2200 1200 1200  5.5200000 5.5300000}%
\special{ar 3400 2200 1200 1200  5.5600000 5.5700000}%
\special{ar 3400 2200 1200 1200  5.6000000 5.6100000}%
\special{ar 3400 2200 1200 1200  5.6400000 5.6500000}%
\special{ar 3400 2200 1200 1200  5.6800000 5.6900000}%
\special{ar 3400 2200 1200 1200  5.7200000 5.7300000}%
\special{ar 3400 2200 1200 1200  5.7600000 5.7700000}%
\special{ar 3400 2200 1200 1200  5.8000000 5.8100000}%
\special{ar 3400 2200 1200 1200  5.8400000 5.8500000}%
\special{ar 3400 2200 1200 1200  5.8800000 5.8900000}%
\special{ar 3400 2200 1200 1200  5.9200000 5.9300000}%
\special{ar 3400 2200 1200 1200  5.9600000 5.9700000}%
\special{ar 3400 2200 1200 1200  6.0000000 6.0100000}%
\special{ar 3400 2200 1200 1200  6.0400000 6.0500000}%
\special{ar 3400 2200 1200 1200  6.0800000 6.0900000}%
\special{ar 3400 2200 1200 1200  6.1200000 6.1300000}%
\special{ar 3400 2200 1200 1200  6.1600000 6.1700000}%
\special{ar 3400 2200 1200 1200  6.2000000 6.2100000}%
\special{ar 3400 2200 1200 1200  6.2400000 6.2500000}%
\special{ar 3400 2200 1200 1200  6.2800000 6.2832853}%
%
\special{pn 8}%
\special{ar 2800 2200 160 680  1.6472752 4.6169857}%
%
\special{pn 8}%
\special{ar 4370 1630 1820 220  1.4694959 1.4812606}%
\special{ar 4370 1630 1820 220  1.5165547 1.5283194}%
\special{ar 4370 1630 1820 220  1.5636136 1.5753783}%
\special{ar 4370 1630 1820 220  1.6106724 1.6224371}%
\special{ar 4370 1630 1820 220  1.6577312 1.6694959}%
\special{ar 4370 1630 1820 220  1.7047900 1.7165547}%
\special{ar 4370 1630 1820 220  1.7518489 1.7636136}%
\special{ar 4370 1630 1820 220  1.7989077 1.8106724}%
\special{ar 4370 1630 1820 220  1.8459665 1.8577312}%
\special{ar 4370 1630 1820 220  1.8930253 1.9047900}%
\special{ar 4370 1630 1820 220  1.9400841 1.9518489}%
\special{ar 4370 1630 1820 220  1.9871430 1.9989077}%
\special{ar 4370 1630 1820 220  2.0342018 2.0459665}%
\special{ar 4370 1630 1820 220  2.0812606 2.0930253}%
\special{ar 4370 1630 1820 220  2.1283194 2.1400841}%
\special{ar 4370 1630 1820 220  2.1753783 2.1871430}%
\special{ar 4370 1630 1820 220  2.2224371 2.2342018}%
\special{ar 4370 1630 1820 220  2.2694959 2.2812606}%
\special{ar 4370 1630 1820 220  2.3165547 2.3283194}%
\special{ar 4370 1630 1820 220  2.3636136 2.3753783}%
\special{ar 4370 1630 1820 220  2.4106724 2.4224371}%
\special{ar 4370 1630 1820 220  2.4577312 2.4694959}%
\special{ar 4370 1630 1820 220  2.5047900 2.5165547}%
\special{ar 4370 1630 1820 220  2.5518489 2.5636136}%
\special{ar 4370 1630 1820 220  2.5989077 2.6106724}%
\special{ar 4370 1630 1820 220  2.6459665 2.6577312}%
\special{ar 4370 1630 1820 220  2.6930253 2.7047900}%
%
\special{pn 8}%
\special{ar 4400 3000 2500 500  3.9507958 3.9587958}%
\special{ar 4400 3000 2500 500  3.9827958 3.9907958}%
\special{ar 4400 3000 2500 500  4.0147958 4.0227958}%
\special{ar 4400 3000 2500 500  4.0467958 4.0547958}%
\special{ar 4400 3000 2500 500  4.0787958 4.0867958}%
\special{ar 4400 3000 2500 500  4.1107958 4.1187958}%
\special{ar 4400 3000 2500 500  4.1427958 4.1507958}%
\special{ar 4400 3000 2500 500  4.1747958 4.1827958}%
\special{ar 4400 3000 2500 500  4.2067958 4.2147958}%
\special{ar 4400 3000 2500 500  4.2387958 4.2467958}%
\special{ar 4400 3000 2500 500  4.2707958 4.2787958}%
\special{ar 4400 3000 2500 500  4.3027958 4.3107958}%
\special{ar 4400 3000 2500 500  4.3347958 4.3427958}%
\special{ar 4400 3000 2500 500  4.3667958 4.3747958}%
\special{ar 4400 3000 2500 500  4.3987958 4.4067958}%
\special{ar 4400 3000 2500 500  4.4307958 4.4387958}%
\special{ar 4400 3000 2500 500  4.4627958 4.4707958}%
\special{ar 4400 3000 2500 500  4.4947958 4.5027958}%
\special{ar 4400 3000 2500 500  4.5267958 4.5347958}%
\special{ar 4400 3000 2500 500  4.5587958 4.5667958}%
\special{ar 4400 3000 2500 500  4.5907958 4.5987958}%
\special{ar 4400 3000 2500 500  4.6227958 4.6307958}%
\special{ar 4400 3000 2500 500  4.6547958 4.6627958}%
\special{ar 4400 3000 2500 500  4.6867958 4.6947958}%
\special{ar 4400 3000 2500 500  4.7187958 4.7267958}%
\special{ar 4400 3000 2500 500  4.7507958 4.7587958}%
%
\special{pn 8}%
\special{pa 2290 3050}%
\special{pa 2700 2740}%
\special{dt 0.045}%
\special{sh 1}%
\special{pa 2700 2740}%
\special{pa 2636 2764}%
\special{pa 2658 2772}%
\special{pa 2660 2796}%
\special{pa 2700 2740}%
\special{fp}%
%
\special{pn 8}%
\special{ar 3290 2200 1310 600  4.7123890 6.2831853}%
%
\special{pn 8}%
\special{ar 3290 2200 1310 600  6.2831853 6.2831853}%
\special{ar 3290 2200 1310 600  0.0000000 1.5707963}%
%
\special{pn 8}%
\special{ar 3290 2200 650 600  1.5707963 4.7123890}%
\end{picture}%
}

\vspace{0.3truecm}

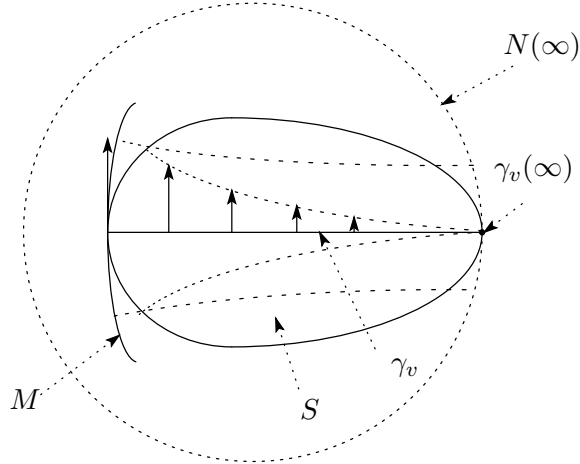
\centerline{{\bf Fig. 2.}}

\section{Proof of Theorem A} 
In this section, we shall prove Theorem A stated in Introduction.  
First we prepare the following lemma.  

\vspace{0.5truecm}

\noindent
{\bf Lemma 4.1.} {\sl Let $M$ be a curvature-adapted submanifold with section 
in a symmetric space $N=G/K$ of non-compact type.  Then the following 
conditions {\rm (i)} and {\rm (ii)} are equivalent:

{\rm(i)} $M$ has no focal point of non-Euclidean type on $N(\infty)$,

{\rm (ii)} for each unit normal vector $v$ of $M$ and each 
$\mu\in{\rm Spec}\,R(\cdot,v)v\setminus\{0\}$, $\pm\sqrt{-\mu}$ 
are not eigenvalues of $A_v\vert_{{\rm Ker}(R(\cdot,v)v-\mu\,{\rm id})}$, 
where $R$ is the curvature tensor of $G/K$ and $A$ is the shape tensor 
of $M$.}

\vspace{0.5truecm}

\noindent
{\it Proof.} First we note that the condition (i) is equivalent to 
the following condition:

\vspace{0.2truecm}

(i$'$) $M$ has no strongly focal point of non-Euclidean type on $N(\infty)$.  

\vspace{0.2truecm}

\noindent
In fact, this fact follows from Proposition 3.1 because $M$ has sections.  
Fix a unit normal vector $v$ of $M$ at any $x=gK\in M$.  
Since $M$ is curvature-adapted, we have 
$$T_xM=\displaystyle{\mathop{\oplus}_{\lambda\in{\rm Spec}\,A_v}
\mathop{\oplus}_{\mu\in{\rm Spec}\,R(\cdot,v)v}\left(
{\rm Ker}(R(\cdot,v)v-\mu\,{\rm id})\cap
{\rm Ker}(A_v-\lambda\,{\rm id})\right),} \leqno{(4.1)}$$
where ${\rm Spec}(\cdot)$ is the spectrum of $(\cdot)$.  
A strongly $M$-Jacobi field 
$Y$ along $\gamma_v$ with $Y(0)\in {\rm Ker}(R(\cdot,v)v-\mu\,{\rm id})\cap
{\rm Ker}(A_v-\lambda\,{\rm id})$ ($\mu\not=0$) is described as 
$$\begin{array}{l}
\displaystyle{Y(t)=P_{\gamma_v\vert_{[0,t]}}\left((D^{co}_{tv}-D^{si}_{tv}
\circ A_{tv})(Y(0))\right)}\\
\hspace{1.2truecm}\displaystyle{=\left(\cosh(t\sqrt{-\mu})
-\frac{\lambda\sinh(t\sqrt{-\mu})}{\sqrt{-\mu}}\right)
P_{\gamma_v\vert_{[0,t]}}(Y(0)).}
\end{array}$$
If $\lambda=\pm\sqrt{-\mu}$, then we have 
$\vert\vert Y(t)\vert\vert=\vert\vert Y(0)\vert\vert 
e^{\pm t\sqrt{-\mu}}$.  Hence we have 
$\lim\limits_{t\to\infty}\frac{\vert\vert Y(t)\vert\vert}{t}=0$ or 
$\lim\limits_{t\to-\infty}\frac{\vert\vert Y(t)\vert\vert}{t}=0$.  
Also, from $Y(0)\in{\rm Ker}(R(\cdot,v)v-\mu\,{\rm id})$ and 
$\mu\not=0$, we have ${\rm Sec}(v,Y(0))<0$.  Hence, either 
$\gamma_v(\infty)$ or $\gamma_{-v}(\infty)$ is a strongly focal point of 
non-Euclidean type on $N(\infty)$ of $M$.  Thus 
${\rm not}\,{\rm (ii)}\,\,\Rightarrow\,\,{\rm not}\,{\rm (i)}$, that is, 
${\rm (i)}\,\,\Rightarrow\,\,{\rm (ii)}$ is shown.  
Assume that (ii) holds.  Take an arbitrary $X(\not=0)\,\in T_xM$ with 
${\rm Sec}(v,X)<0$.  Let $S:=\{(\lambda,\mu)\in{\rm Spec}\,A_v\times
{\rm Spec}\,R(\cdot,v)v\,\vert\,{\rm Ker}(R(\cdot,v)v-\mu\,{\rm id})\cap
{\rm Ker}(A_v-\lambda\,{\rm id})\not=\{0\}\}$ and 
$S_0:=\{(\lambda,\mu)\in S\,\vert\,\mu=0\}$.  
Let $X=\sum\limits_{(\lambda,\mu)\in S}X_{\lambda,\mu}$, where 
$X_{\lambda,\mu}\in{\rm Ker}(R(\cdot,v)v-\mu\,{\rm id})\cap
{\rm Ker}(A_v-\lambda\,{\rm id})$.  The strongly $M$-Jacobi field $Y$ along 
$\gamma_v$ with $Y(0)=X$ is described as 
$$\begin{array}{l}
\displaystyle{Y(t)=P_{\gamma_v\vert_{[0,t]}}\left((D^{co}_{tv}-D^{si}_{tv}
\circ A_{tv})(X)\right)}\\
\hspace{1truecm}\displaystyle{=\sum_{(\lambda,\mu)\in S\setminus S_0}
\left(\cosh(t\sqrt{-\mu})-\frac{\lambda\sinh(t\sqrt{-\mu})}
{\sqrt{-\mu}}\right)P_{\gamma_v\vert_{[0,t]}}(X_{\lambda,\mu})}\\
\hspace{1.4truecm}\displaystyle{+\sum_{(\lambda,\mu)\in S_0}
(1-t\lambda)P_{\gamma_v\vert_{[0,t]}}(X_{\lambda,\mu})}.
\end{array}$$
Hence we have 
$$\begin{array}{l}
\displaystyle{\vert\vert Y(t)\vert\vert^2=
\sum_{(\lambda,\mu)\in S\setminus S_0}
\left(\cosh(t\sqrt{-\mu})
-\frac{\lambda\sinh(t\sqrt{-\mu})}
{\sqrt{-\mu}}\right)^2\vert\vert X_{\lambda,\mu}\vert\vert^2}\\
\hspace{2truecm}\displaystyle{+\sum_{(\lambda,\mu)\in S_0}
(1-t\lambda)^2\vert\vert X_{\lambda,\mu}\vert\vert^2.}
\end{array}$$
Since ${\rm Sec}(v,X)<0$, there exists $(\lambda_0,\mu_0)\in S\setminus S_0$ 
with $X_{\lambda_0,\mu_0}\not=0$.  Then we have 
$$\begin{array}{l}
\displaystyle{\frac{\vert\vert Y(t)\vert\vert^2}{t^2}\geq
\frac{1}{t^2}\left(\cosh(t\sqrt{-\mu_0})
-\frac{\lambda_0\sinh(t\sqrt{-\mu_0})}
{\sqrt{-\mu_0}}\right)^2\vert\vert X_{\lambda_0,\mu_0}\vert\vert^2}\\
\hspace{1.1truecm}\displaystyle{=\frac{1}{2t^2}
\left((1-\frac{\lambda_0}{\sqrt{-\mu_0}})
e^{t\sqrt{-\mu_0}}+(1+\frac{\lambda_0}{\sqrt{-\mu_0}})
e^{-t\sqrt{-\mu_0}}\right)^2
\vert\vert X_{\lambda_0,\mu_0}\vert\vert^2}.
\end{array}$$
Hence, since $\lambda_0\not=\pm\sqrt{-\mu_0}$ by the assumption, 
we have $\lim\limits_{t\to\infty}\dfrac{\vert\vert Y(t)\vert\vert}{t}=\infty$.  From the arbitrariness of $X\in T_xM$, it follows that $\gamma_v(\infty)$ is 
not a focal point of non-Euclidean type (on $N(\infty)$) of $M$.  
Furthermore, from the arbitrarinesses of $v$ and $x$, we see that $M$ has no 
focal point of non-Euclidean type on $N(\infty)$.  Thus 
${\rm (ii)}\,\Rightarrow\,{\rm (i)}$ is shown.  \hspace{1.5truecm}q.e.d.

\vspace{0.5truecm}

\noindent
{\it Remark 4.1.} By imitating the proof of this lemma, it is shown that the 
following conditions (i$'$) and (ii$'$) are equivalent:

(i$'$) $M$ has no focal point on $N(\infty)$, 

(ii$'$) for each unit normal vector $v$ of $M$ and each 
$\mu\in{\rm Spec}\,R(\cdot,v)v$, $\pm\sqrt{-\mu}$ are not eigenvalues of 
$A_v\vert_{{\rm Ker}(R(\cdot,v)v-\mu I)}$.  

\vspace{0.5truecm}

By using this lemma, the statement of Theorem A is proved.  

\vspace{0.5truecm}

\noindent
{\it Proof of Theorem A.} According to the proof of the statement (ii) of 
Theorem 1 in [Koi2], $M$ is proper complex equifocal if and only if the 
following condition ($\ast$) holds:

\vspace{0.2truecm}

($\ast$) for each unit normal vector $v$ of $M^{\bf c}$ and each 
$\mu\in{\rm Spec}_JR^{\bf c}(\cdot,v)v\setminus\{0\}$, 
$\sqrt{\mu}$ ($2$-values) are not $J$-eigenvalues of 
$A^{\bf c}_v\vert_{{\rm Ker}(R^{\bf c}(\cdot,v)v-\mu\,{\rm id})}$, 
where $A^{\bf c}$ is the shape tensor of $M^{\bf c}$, $R^{\bf c}$ is 
the curvature tensor of $G^{\bf c}/K^{\bf c}$, $J$ is the complex 
structures of $M^{\bf c}$ and 
${\rm Spec}_J(\cdot)$ is the $J$-spectrum of $(\cdot)$.  

\vspace{0.2truecm}

\noindent
It is easy to show that this condition ($\ast$) is equivalent to the condition 
(ii) of Lemma 4.1.  Hence the statement of Theorem A follows from Lemma 4.1.  
\hspace{3.5truecm} q.e.d.

\section{Proofs of Theorems B$\sim$E} 
In this section, we shall prove Theorems B$\sim$E.  
For its purpose, we prepare a lemma.  
Let $M$ be a curvature-adapted and proper complex equifocal 
$C^{\omega}$-submanifold in a 
symmetric space $G/K$ of non-compact type, 
where we may assume that $eK\in M$ ($e\,:\,$ the identity element of $G$) 
by operating an element of $G$ to $M$ if necessary 
and hence the constant path $\hat 0$ at the zero element $0$ of 
$\mathfrak g$ is contained in $\widetilde{M}:=
(\pi\circ\phi)^{-1}(M)$.  
Denote by $\widetilde M_0$ the component of $\widetilde M$ containing 
$\hat 0$.  Fix a unit normal vector 
$v$ of $M$ at $eK$.  Set $\mathfrak p:=T_{eK}(G/K)$ and 
$\mathfrak b:=T^{\perp}_{eK}M$.  Let $\mathfrak p=\mathfrak a+
\sum\limits_{\alpha\in\triangle_+}\mathfrak p_{\alpha}$ be the root space 
decomposition with respect to a lexicographically ordered maximal abelian 
subspace $\mathfrak a$ containing $\mathfrak b$.  
Let $\overline{\triangle}:=\{\alpha\vert_{\mathfrak b}\,\vert\,\alpha\in
\triangle\,\,{\rm s.t.}\,\,\alpha\vert_{\mathfrak b}\not=0\}$ and 
$\mathfrak p=\mathfrak z_{\mathfrak p}(\mathfrak b)
+\sum_{\beta\in\overline{\triangle}_+}
\mathfrak p_{\beta}$ be the root space decomposition with respect to 
$\mathfrak b$, where $\mathfrak z_{\mathfrak p}(\mathfrak b)$ is 
the centralizer of $\mathfrak b$ in $\mathfrak p$.  
For convenience, we denote $\mathfrak z_{\mathfrak p}(\mathfrak b)$ 
by $\mathfrak p_0$.  
Then we have 
$\mathfrak p_{\beta}=\sum_{\alpha\in\triangle_+\,\,{\rm s.t.}\,\,
\alpha\vert_{\mathfrak b}=\pm\beta}\mathfrak p_{\alpha}$ 
($\beta\in\overline{\triangle}_+$) and 
$\mathfrak p_0=\mathfrak a+\sum_{\alpha\in\triangle_+\,\,{\rm s.t.}\,\,
\alpha\vert_{\mathfrak b}=0}\mathfrak p_{\alpha}$.  
Let $v^L$ be the horizontal lift of $v$ to $\hat 0$.  Denote by $A$ 
(resp. $\widetilde A$) the shape tensor of $M$ (resp. $\widetilde{M}_0$).  
According to Theorem 5.9 of [Koi1], we have the following fact.  

\vspace{0.5truecm}

\noindent
{\bf Lemma 5.1.} {\sl If the spectrum of $A_v$ is equal to 
$\{\lambda_1,\cdots,\lambda_g\}$, then the spectrum of 
$\widetilde A_{v^L}^{\bf c}$ is given by 
\begin{align*}
&\{0\}\cup\{\lambda_i\,\vert\,i\in I_0\}\\
&\cup\left(\mathop{\cup}_{\mu\in{\rm Spec}\,R(\cdot,v)v\setminus\{0\}}
\{\frac{\sqrt{-\mu}}{{\rm arctanh}\frac{\sqrt{-\mu}}{\lambda_i}+j\pi
\sqrt{-1}}\,\vert\,
i\in I_{\mu}^+,\,\,j\in {\bf Z}\}\right)\\
&\cup\left(\mathop{\cup}_{\mu\in{\rm Spec}\,R(\cdot,v)v\setminus\{0\}}
\{\frac{\sqrt{-\mu}}{{\rm arctanh}\frac{\lambda_i}{\sqrt{-\mu}}+
(j+\frac12)\pi\sqrt{-1}}\,\vert\,
i\in I_{\mu}^-,\,\,j\in{\bf Z}\}\right),
\end{align*}
where $I_0=\{i\,\vert\,{\rm Ker}\,R(\cdot,v)v\cap{\rm Ker}(A_v-\lambda_i\,
{\rm id})\not=\{0\}\}$, 
$I_{\mu}^+:=\{i\in I_{\mu}\,\vert\,\vert\lambda_i\vert\,>\,\sqrt{-\mu}\}$ 
and $I_{\mu}^-:=\{i\in I_{\mu}\,\vert\,\vert\lambda_i\vert\,<\,\sqrt{-\mu}\}$ 
as $I_{\mu}:=\{i\,\vert\,{\rm Ker}(R(\cdot,v)v-\mu\,{\rm id})\cap
{\rm Ker}(A_v-\lambda_i\,{\rm id})\not=\{0\}\}$.}

\vspace{0.5truecm}

Now we shall prove Theorems B and C in terms of this lemma and Lemma 2.1.  

\vspace{0.5truecm}

\noindent
{\it Proof of Theorems B and C.} 
Let $m_A:=\displaystyle{\mathop{\max}_{v\in\mathfrak b\setminus\{0\}}
\sharp{\rm Spec}\,A_v}$ and 
$m_R:=\displaystyle{\mathop{\max}_{v\in\mathfrak b\setminus\{0\}}
\sharp{\rm Spec}\,R(\cdot,v)v}$.  
Let $U:=\{v\in\mathfrak b\setminus\{0\}\,\vert\,
\sharp{\rm Spec}A_v=m_A,\,\,\sharp{\rm Spec}\,R(\cdot,v)v=m_R\}$, 
which is an open dense subset of $\mathfrak b\setminus\{0\}$.  
Note that ${\rm Spec}\,R(\cdot,v)v=\{-\beta(v)^2\,\vert\,\beta\in
\overline{\triangle}_+\}$ and, if $v\in U$, then $\beta(v)^2$'s 
($\beta\in\overline{\triangle}_+$) are mutually distinct (i.e., 
$m_R=\sharp\overline{\triangle}_+$).  
Let ${\rm Spec}A_v=\{\lambda^v_1,\cdots,\lambda^v_{m_A}\}$ 
($\lambda^v_1>\cdots>\lambda^v_{m_A}$) ($v\in U$).  
Then it follows from Lemma 5.1 that 
$$\begin{array}{l}
\hspace{0.5truecm}\displaystyle{{\rm Spec}\widetilde A_{v^L}^{\bf c}}\\
\displaystyle{=\{0\}\cup\{\lambda^v_i\,\vert\,i\in I_0^v\}}\\
\hspace{0.5truecm}\displaystyle{\cup\left(
\mathop{\cup}_{\beta\in\overline{\triangle}_+}
\{\frac{\beta(v)}{{\rm arctanh}\frac{\beta(v)}{\lambda^v_i}
+j\pi\sqrt{-1}}\,\vert\,i\in (I_{\beta}^v)^+,\,\,j\in {\bf Z}\}\right)}\\
\hspace{0.5truecm}\displaystyle{\cup\left(
\mathop{\cup}_{\beta\in\overline{\triangle}_+}
\{\frac{\beta(v)}{{\rm arctanh}\frac{\lambda^v_i}{\beta(v)}+
(j+\frac12)\pi\sqrt{-1}}\,\vert\,
i\in (I_{\beta}^v)^-,\,\,j\in{\bf Z}\}\right)}
\end{array}\leqno{(5.1)}$$
for any $v\in U$, where 
$I_0^v:=\{i\,\vert\,\mathfrak p_0\cap{\rm Ker}(A_v-\lambda_i^v{\rm id})\not=
\{0\}\}$, 
$(I_{\beta}^v)^+:=\{i\in I_{\beta}^v\,\vert\,
\vert\lambda^v_i\vert\,>\,\vert\beta(v)\vert\}$ and 
$(I_{\alpha}^v)^-:=\{i\in I_{\beta}^v\,\vert\,\vert\lambda^v_i\vert\,<\,
\vert\beta(v)\vert\}$ as 
$I_{\beta}^v:=\{i\,\vert\,\mathfrak p_{\beta}\cap{\rm Ker}(A_v-\lambda^v_i
{\rm id})\not=\{0\}\}$.  
Let $F$ be the sum of all complex focal hyperplanes of 
$(\widetilde M_0,\widehat 0)$.  From $(5.1)$, the set 
$$\mathop{\cup}_{v\in U}
\left(
\begin{array}{l}
\displaystyle{\{\frac{1}{\lambda^v_i}v^L\,\vert\,i\in I_0^v\,\,{\rm s.t.}\,\,
\lambda_i^v\not=0\}\cup}\\
\displaystyle{\left(\mathop{\cup}_{\beta\in\overline{\triangle}_+}
\{\frac{{\rm arctanh}\frac{\beta(v)}{\lambda_i^v}+j\pi\sqrt{-1}}{\beta(v)}v^L
\,\vert\,\,i\in(I_{\beta}^v)^+,j\in {\bf Z}\}\right)\cup}\\
\displaystyle{\left(
\mathop{\cup}_{\beta\in\overline{\triangle}_+}
\{\frac{{\rm arctanh}\frac{\lambda^v_i}{\beta(v)}+(j+\frac12)\pi\sqrt{-1}}
{\beta(v)}v^L\,\vert\,i\in(I_{\beta}^v)^-,j\in{\bf Z}\}\right)}
\end{array}\right)
\leqno{(5.2)}$$
is contained in $F$.  Fix $v_0\in U$.  Since the set $(5.2)$ 
is contained in $F$ and $F$ consists of infinitely many complex 
hyperplanes of $(T^{\perp}_{\hat0}\widetilde M_0)^{\bf c}$, 
it is shown by delicate argument that there exist the complex linear 
functions $\phi_i$ ($i\in I_0^{v_0}\,\,{\rm s.t.}\,\,\lambda_i^{v_0}\not=0$), 
$\phi^1_{\beta,i,j}$ ($\beta\in\overline{\triangle}_+,\,i\in(I^{v_0}_{\beta})^+,\,j\in{\bf Z}$) and 
$\phi^2_{\beta,i,j}$ ($\beta\in\overline{\triangle}_+,\,
i\in(I^{v_0}_{\beta})^-,\,j\in{\bf Z}$) on 
$(T^{\perp}_{\hat0}\widetilde M_0)^{\bf c}(={\mathfrak b}^{\bf c})$ 
satisfying 
$\phi_i(v)=\lambda^v_i\,\,(v\in U')$, 
$\displaystyle{\phi^1_{\beta,i,j}(v)=
\frac{\beta(v)}{{\rm arctanh}\frac{\beta(v)}{\lambda_i^v}+j\pi\sqrt{-1}}
\,(v\in U')}$ and 
$\displaystyle{
\phi^2_{\beta,i,j}(v)=}$\newline
$\displaystyle{
\frac{\beta(v)}{{\rm arctanh}
\frac{\lambda_i^v}{\beta(v)}+(j+\frac12)\pi\sqrt{-1}}}$ $(v\in U')$, 
respectively, where $U'$ is a sufficiently small neighborhood of $v_0$ 
in $U$.  
Since $\phi^k_{\beta,i,j}(v)=\frac{\beta(v)\phi^k_{\beta,i,0}(v)}
{\beta(v)+j\pi\phi^k_{\beta,i,0}(v)\sqrt{-1}}$ for all $v\in U'$ and all 
$j\in{\bf Z}$ and $\phi^k_{\beta,i,j}$'s are complex linear, we see that 
$\frac{\beta(v)}{\phi^k_{\beta,i,j}(v)}$ is independent of the choice of 
$v\in U'$, where $\beta\in\overline{\triangle}_+$ and 
$(k,i)\in(\{1\}\times(I_{\beta}^{v_0})^+)\cup
(\{2\}\times(I_{\beta}^{v_0})^-)$.  That is, 
$\displaystyle{\frac{\beta(v)}{\lambda_i^v}}$ ($i\in(I_{\beta}^{v_0})^+\,
(\beta\in\overline{\triangle}_+)$) and $\displaystyle{\frac{\lambda_i^v}
{\beta(v)}}$ ($i\in(I_{\beta}^{v_0})^-\,(\beta\in\overline{\triangle}_+$)) 
are independent of the choices of $v\in U'$.  Set 
$\displaystyle{c^+_{\beta,i}:=\frac{\beta(v_0)}{\lambda_i^{v_0}}}$ 
$(i\in(I_{\beta}^{v_0})^+\,(\beta\in\overline{\triangle}_+))$ 
and $\displaystyle{c^-_{\beta,i}:=\frac{\lambda_i^{v_0}}{\beta(v_0)}}$ 
$(i\in(I_{\beta}^{v_0})^-\,(\beta\in\overline{\triangle}_+))$.  
Hence we have 
$\displaystyle{\phi^1_{\beta,i,j}=\frac{\beta^{\bf c}
\vert_{\mathfrak b^{\bf c}}}
{{\rm arctanh}c^+_{\beta,i}+j\pi\sqrt{-1}}}$ and 
$\displaystyle{\phi^2_{\beta,i,j}
=\frac{\beta^{\bf c}\vert_{\mathfrak b^{\bf c}}}
{{\rm arctanh}c^-_{\beta,i}+(j+\frac12)\pi\sqrt{-1}}}$.  
Clearly we have 
$$\begin{array}{l}
\displaystyle{F=\left(\mathop{\cup}_{i\in I_0^{v_0}\,{\rm s.t.}\,
\lambda_i^{v_0}\not=0}\phi_i^{-1}(1)\right)}\\
\hspace{1truecm}\displaystyle{\cup\left(
\mathop{\cup}_{\beta\in\overline{\triangle}_+}
\mathop{\cup}_{(i,j)\in(I^{v_0}_{\beta})^+\times{\bf Z}}
(\phi^1_{\beta,i,j})^{-1}(1)\right)}\\
\hspace{1truecm}\displaystyle{\cup\left(
\mathop{\cup}_{\beta\in\overline{\triangle}_+}
\mathop{\cup}_{(i,j)\in(I^{v_0}_{\beta})^-\times{\bf Z}}
(\phi^2_{\beta,i,j})^{-1}(1)\right)}\\
\hspace{0.4truecm}
\displaystyle{=\left(
\mathop{\cup}_{i\in I_0^{v_0}\,{\rm s.t.}\,
\lambda_i^{v_0}\not=0}\phi_i^{-1}(1)\right)}\\
\hspace{1truecm}\displaystyle{\cup\left(
\mathop{\cup}_{\beta\in\overline{\triangle}_+}
\mathop{\cup}_{(i,j)\in(I^{v_0}_{\beta})^+\times{\bf Z}}
(\beta^{\bf c})^{-1}({\rm arctanh}c_{\beta,i}^++j\pi\sqrt{-1})\right)}\\
\hspace{1truecm}\displaystyle{\cup\left(
\mathop{\cup}_{\beta\in\overline{\triangle}_+}
\mathop{\cup}_{(i,j)\in(I^{v_0}_{\beta})^-\times{\bf Z}}
(\beta^{\bf c})^{-1}({\rm arctanh}c_{\beta,i}^-+(j+\frac12)\pi\sqrt{-1})
\right).}
\end{array}
\leqno{(5.3)}$$
Also, we can show that the complex reflection 
group generated by the complex reflections of order two with respect to 
$(\beta^{\bf c})^{-1}(0)$'s ($\beta\in\overline{\triangle}_+$) is of rank $r$, where $r:={\rm codim}\,M$.  
The group $W_{\widetilde M_0}$ is generated by the complex reflections of 
order two with respect to the complex hyperplanes in $(5.3)$ constructing 
$F$.  This group is isomorphic to the complex Coxeter group $W_M$ associated 
with $M$ and hence it is discrete and, according to Lemma 3.5 of [Koi4], $F$ 
is $W_{\widetilde M_0}$-invariant.  Therefore, we have 
$$\begin{array}{l}
\displaystyle{F=\left(
\mathop{\cup}_{\beta\in\overline{\triangle}_+}
\mathop{\cup}_{(i,j)\in(I^{v_0}_{\beta})^+\times{\bf Z}}
(\beta^{\bf c})^{-1}({\rm arctanh}c_{\beta,i}^++j\pi\sqrt{-1})\right)}\\
\hspace{1truecm}\displaystyle{\cup\left(
\mathop{\cup}_{\beta\in\overline{\triangle}_+}
\mathop{\cup}_{(i,j)\in(I^{v_0}_{\beta})^-\times{\bf Z}}
(\beta^{\bf c})^{-1}({\rm arctanh}c_{\beta,i}^-+(j+\frac12)\pi\sqrt{-1})
\right),}
\end{array}
\leqno{(5.4)}$$
where we note that $\{i\in I_0^{v_0}\,\vert\,\lambda_i^{v_0}\not=0\}$ 
is not necessarily empty set.  
Denote by ${\rm proj}_{\bf R}$ the natural projection of 
$\mathfrak b^{\bf c}$ onto $\mathfrak b$ and set 
$F_{\bf R}:={\rm proj}_{\bf R}(F)$.  Then we have 
$$\begin{array}{l}
\displaystyle{F_{\bf R}=\left(
\mathop{\cup}_{\beta\in\overline{\triangle}_+}
\mathop{\cup}_{i\in(I^{v_0}_{\beta})^+}
\beta^{-1}({\rm arctanh}c_{\beta,i}^+)\right)}\\
\hspace{1truecm}\displaystyle{\cup\left(
\mathop{\cup}_{\beta\in\overline{\triangle}_+}
\mathop{\cup}_{i\in(I^{v_0}_{\beta})^-}
\beta^{-1}({\rm arctanh}c_{\beta,i}^-)\right).}
\end{array}
\leqno{(5.5)}$$
Let $W'_{\widetilde M_0}$ be the group generated by the reflections 
(in $\mathfrak b$) with respect to the hyperplanes 
constructing $F_{\bf R}$.  Since $F$ is 
$W_{\widetilde M_0}$-invariant, $F_{\bf R}$ is 
$W'_{\widetilde M_0}$-invariant.  Therefore, since $F_{\bf R}$ consists of 
finite pieces of (real) hyperplanes (in $\mathfrak b$), the intersection 
of all the hyperplanes constructing $F_{\bf R}$ is non-empty.  
Take an element $Z$ of the intersection.  
Then complex hyperplanes in $(5.4)$ constructing $F$ are rewritten as 
$$\begin{array}{l}
\displaystyle{(\beta^{\bf c})^{-1}
({\rm arctanh}c_{\beta,i}^++j\pi\sqrt{-1})=Z+(\beta^{\bf c})^{-1}
(j\pi\sqrt{-1}),}\\
\displaystyle{(\beta^{\bf c})^{-1}
({\rm arctanh}c_{\beta,i}^-+(j+\frac12)\pi\sqrt{-1})=
Z+(\beta^{\bf c})^{-1}((j+\frac12)\pi\sqrt{-1}),}
\end{array}
\leqno{(5.6)}$$
respectively.  Hence we see that $W_{\widetilde M_0}$ is isomorphic to 
the group generated by the (real) reflections with 
respect to the hyperplanes 
${\widehat{\beta}}^{-1}(j\pi)$'s ($\beta\in\overline{\triangle}^+_+,\,
j\in{\bf Z}$) 
and ${\widehat{\beta}}^{-1}((j+\frac12)\pi)$'s 
($\beta\in\overline{\triangle}^-_+,\,j\in{\bf Z}$) in 
$\sqrt{-1}\mathfrak b$, 
where $\widehat{\beta}:=-\sqrt{-1}\beta^{\bf c}\vert_{\sqrt{-1}\mathfrak b}$ 
and $\overline{\triangle}^{\pm}_+:=
\{\beta\in\overline{\triangle}_+\,\vert\,
(I^{v_0}_{\beta})^{\pm}\not=\emptyset\}$.  
Thus $W_{\widetilde M_0}$ is isomorphic to the affine transformation group 
associated with $\overline{\triangle}$.  Hence, since $F$ is 
$W_{\widetilde M_0}$-invariant, we see that $\overline{\triangle}$ 
is a weakly root system.  
This completes the proof of Theorem B.  
According to $(5.6)$, for each fixed $\beta\in\overline{\triangle}_+$, 
$c_{\beta,i}^+$'s ($i\in(I_{\beta}^{v_0})^+$) coincide and so are 
$c_{\beta,i}^-$'s ($i\in(I_{\beta}^{v_0})^-$) also.  
In particular, we have $\sharp(I_{\beta}^{v_0})^+\leq1$ and 
$\sharp(I_{\beta}^{v_0})^-\leq 1$.  
This fact implies that 
$\sharp{\rm Spec}\,A_{v_0}$ is evaluated from above as in the statement of 
Theorem C.  From $v_0\in U$ and the definition of $U$, 
it follows that $\sharp{\rm Spec}\,A_v\leq\sharp{\rm Spec}\,A_{v_0}$ 
for any normal vector $v$.  
Therefore the statement of Theorem C follows.  \hspace{3.8truecm} q.e.d.

\vspace{0.5truecm}

\noindent
{\it Remark 5.1.} In the case where $M$ is curvature-adapted equifocal 
submanifold in a symmetric space $G/K$ of compact type, we have 
$$\begin{array}{l}
\hspace{0.5truecm}\displaystyle{{\rm Spec}\widetilde A_{v^L}^{\bf c}=
{\rm Spec}\widetilde A_{v^L}}\\
\displaystyle{=\{0\}\cup\{\lambda^v_i\,\vert\,i\in I_0^v\}}\\
\hspace{0.5truecm}\displaystyle{\cup\left(
\mathop{\cup}_{\beta\in\overline{\triangle}_+}
\{\frac{\beta(v)}{{\rm arctanh}\frac{\beta(v)}{\lambda^v_i}
+j\pi}\,\vert\,i\in I_{\beta}^v,\,\,j\in {\bf Z}\}\right),}
\end{array}$$
where $\widetilde A_{v^L},\,\overline{\triangle}_+, I_0^v$ and $I_{\beta}^v$ 
are as in the above proof.  Also, we have 
$$F=\mathop{\cup}_{\beta\in\overline{\triangle}_+}
\mathop{\cup}_{(i,j)\in I^{v_0}_{\beta}\times{\bf Z}}
(\beta^{\bf c})^{-1}({\rm arctanh}c_{\beta,i}+j\pi)$$
and hence 
$$F_{\bf R}=\mathop{\cup}_{\beta\in\overline{\triangle}_+}
\mathop{\cup}_{(i,j)\in I^{v_0}_{\beta}\times{\bf Z}}
\beta^{-1}({\rm arctanh}c_{\beta,i}+j\pi),\leqno{(5.7)}$$
where 
$F,F_{\bf R}$ and $v_0$ are as in the above proof and 
$c_{\beta,i}:=\frac{\beta(v_0)}{\lambda_i^{v_0}}$.  
Furthermore, it is shown that $F_{\bf R}$ is $W'_{\widetilde M_0}$-invariant, 
where $W'_{\widetilde M_0}$ is as in the above proof.  Note that 
$W'_{\widetilde M_0}$ is the affine Coxeter group asociated with the 
isoparametric submanifold $\widetilde M_0$.  However, it does not follow from 
these facts that, for each fixed $\beta$, $c_{\beta,i}$'s 
($i\in I_{\beta}^{v_0}$) coincide because of the existenceness of the term 
$j\pi$ in the right-hand side of the relation $(5.7)$.  Thus we cannot 
evaluate $\sharp{\rm Spec}\,A_v$ from above 
for curvature-adapted equifocal submanifolds in a symmetric space of 
compact type.  

\vspace{0.5truecm}

Next we shall prove Corollary B.1.  

\vspace{0.5truecm}

\noindent
{\it Proof of Corollary B.1.}  According to Theorem 2 in [Koi4], it follows 
that $M$ is decomposed into the (non-trivial) product of two 
curvature adapted and proper complex equifocal submanifolds if and only if 
the complex Coxeter group associated with $M$ is decomposable.  
Hence the statement of Corollary B.1 follows from Theorem B.  
\hspace{3.7truecm}q.e.d.

\vspace{0.5truecm}

Next we shall prove Corollary B.2 in terms of Corollary B.1.  

\vspace{0.5truecm}

\noindent
{\it Proof of Corollary B.2.} Since ${\rm codim}\,M={\rm rank}\,G/K$, we have 
$\overline{\triangle}=\triangle$, that is, 
$W_{\overline{\triangle}}$ is equal to the Weyl group associated with the 
symmetric space $G/K$.  Hence, since $G/K$ is reducible, 
$W_{\overline{\triangle}}$ is decomposable.  
Therefore, the statement of Corollary B.2 follows from Corollary B.1.  
\hspace{11.5truecm}q.e.d.

\vspace{0.5truecm}

Next we shall prove Theorem D.  

\vspace{0.5truecm}

\noindent
{\it Proof of Theorem D.} Without loss of generality, we may assume 
$x_0=eK$.  According to the proof of Theorem B, the sum $F$ of 
all complex focal hyperplanes of 
$(\widetilde M,\hat0)$ is as in $(5.4)$.  
The intersection of 
$F\,(\subset (T^{\perp}_{\hat0}\widetilde M)^{\bf c}
={\mathfrak b}^{\bf c})$ 
with $\mathfrak b$ is as follows:
$$F\cap \mathfrak b=
\mathop{\cup}_{\beta\in\overline{\triangle}_+}
\mathop{\cup}_{i\in(I_{\beta}^{v_0})^+}\beta^{-1}
({\rm arctanh}c_{\beta,i}^+).
\leqno{(5.8)}$$
Since $\beta^{-1}({\rm arctanh}c_{\beta,i}^+)$ 
($i\in(I_{\beta}^{v_0})^+$ ($\beta\in\overline{\triangle}_+$)) 
are (real) hyperplanes in $\mathfrak b$ through $Z$ in the proof of 
Theorems B and C and 
$\mathfrak b\,(\subset\mathfrak p\subset\mathfrak g)$ is abelian, 
$\exp^{\perp}(\beta^{-1}({\rm arctanh}c_{\beta,i}^+))$ 
($i\in(I_{\beta}^{v_0})^+$ ($\beta\in\overline{\triangle}_+$)) 
are totally geodesic hypersurfaces through $\exp^{\perp}(Z)$ in the section 
$\Sigma:=\exp^{\perp}(\mathfrak b)$.  On the other hand, it is clear that 
$\exp^{\perp}(F\cap\mathfrak b)$ is the focal set of $(M,eK)$.  
Hence, the statement of Theorem D follows.  \hspace{6.4truecm}q.e.d.

\vspace{0.5truecm}

Next we shall prove Theorem E.  

\vspace{0.5truecm}

\noindent
{\it Proof of Theorem E.}  
From $(5.4)$ and $(5.8)$, the statement of Theorem E follows.  
\hspace{0.5truecm}  q.e.d.

\vspace{1truecm}

\centerline{{\bf Appendix 1}}

\vspace{0.5truecm}

In this appendix, we shall first calculate the complex Coxeter group $W_M$ and 
the real Coxeter group $W_{M,{\bf R}}$ associated with a principal orbit $M$ 
of a Hermann type action $H\times G/K\to G/K$ without use of Theorems B 
and E.  Let $\theta$ be the Cartan involution of $G$ with 
$({\rm Fix}\,\theta)_0\subset K\subset{\rm Fix}\,\theta$ and $\sigma$ be 
an involution of $G$ with $({\rm Fix}\,\sigma)_0\subset H\subset{\rm Fix}\,
\sigma$.  Without loss of generality, we may assume that $\sigma\circ\theta
=\theta\circ\sigma$.  Denote by $A$ the shape tensor of $M$.  Then $H(eK)$ is 
a totally geodesic singular orbit of the $H$-action and $M$ is catched as 
a partial tube over $H(eK)$.  Let $L:={\rm Fix}(\sigma\circ\theta)$.  The 
submanifold $\exp^{\perp}(T^{\perp}_{eK}(H(eK)))$ is totally geodesic and it 
is isometric to the symmetric space $L/H\cap K$, where $\exp^{\perp}$ is the 
normal exponential map of $H(eK)$.  
Let $\mathfrak g,\,\mathfrak f$ and $\mathfrak h$ be the Lie algebras of 
$G,\,K$ and $H$.  Denote by the same symbols the involutions of $\mathfrak g$ 
associated with $\theta$ and $\sigma$.  Set $\mathfrak p:={\rm Ker}
(\theta+{\rm id})\,(\subset\mathfrak g)$ and $\mathfrak q:={\rm Ker}
(\sigma+{\rm id})\,(\subset\mathfrak g)$.  
Take $x:=\exp^{\perp}(\xi)=\exp_G(\xi)K
\in M\cap\exp^{\perp}(T^{\perp}_{eK}(H(eK)))$, where $\xi\in\mathfrak p$.  
For simplicity, set $g:=\exp_G(\xi)$.  
Let $\Sigma$ be the section of $M$ through $x$, which pass 
through $eK$.  Let $\mathfrak b:=T_{eK}\Sigma$, $\mathfrak a$ be a maximal 
abelian subspace of $\mathfrak p=T_{eK}(G/K)$ containing $\mathfrak b$, 
$\triangle$ be the root system with respect to $\mathfrak a$ and 
$\mathfrak p=\mathfrak a+\sum\limits_{\alpha\in\triangle_+}
\mathfrak p_{\alpha}$ be the root space decomposition with respect to 
$\mathfrak a$.  
Set $\mathfrak p':=\mathfrak p\cap\mathfrak q (=T^{\perp}_{eK}(H(eK)))$.  
The orthogonal complement ${\mathfrak p'}^{\perp}$ of $\mathfrak p'$ in 
$\mathfrak p$ is equal to $\mathfrak p\cap\mathfrak h$.  
Set $\overline{\triangle}:=\{\alpha\vert_{\mathfrak b}\,
\vert\,\alpha\in\triangle\,\,{\rm s.t.}\,\,\alpha\vert_{\mathfrak b}\not=0\}$, 
$\mathfrak p_{\beta}:=\sum\limits_
{\alpha\in\triangle_+\,\,{\rm s.t.}\,\,\alpha\vert_{\mathfrak b}=\pm\beta}
\mathfrak p_{\alpha}$ for $\beta\in\overline{\triangle}_+$, 
$\overline{\triangle}_+^H:=\{\beta\in\overline{\triangle}_+\,\vert\,
{\mathfrak p'}^{\perp}\cap\mathfrak p_{\beta}\not=\{0\}\}$ and 
$\overline{\triangle}_+^V:=\{\beta\in\overline{\triangle}_+\,\vert\,
\mathfrak p'\cap\mathfrak p_{\beta}\not=\{0\}\}$.  Since both $\mathfrak p'$ 
and ${\mathfrak p'}^{\perp}$ are Lie triple systems of $\mathfrak p$ and 
$\mathfrak b$ is contained in $\mathfrak p'$, we have 
${\mathfrak p'}^{\perp}=\mathfrak z_{{{\mathfrak p}'}^{\perp}}(\mathfrak b)+
\sum\limits_{\beta\in\overline{\triangle}_+^H}
({\mathfrak p'}^{\perp}\cap\mathfrak p_{\beta})$ and 
$\mathfrak p'=\mathfrak b+\sum\limits_{\beta\in\overline{\triangle}_+^V}
(\mathfrak p'\cap\mathfrak p_{\beta})$, where 
$\mathfrak z_{{{\mathfrak p}'}^{\perp}}(\mathfrak b)$ is the centralizer of 
$\mathfrak b$ in ${{\mathfrak p}'}^{\perp}$.  
Take $\eta\in T^{\perp}_xM$.  For each $X\in {\mathfrak p'}^{\perp}\cap
\mathfrak p_{\beta}$ ($\beta\in\overline{\triangle}^H_+$), we can show 
$$A_{\eta}\widetilde X_{\xi}=-\beta(\bar{\eta})\tanh\beta(\xi)
\widetilde X_{\xi}\leqno{({\rm A}.1)}$$
(see the proof of Theorem B of [Koi3]), where $\widetilde X_{\xi}$ is the 
horizontal lift of $X$ to $\xi$ (see Section 3 of [Koi3] about this 
definition) and $\bar{\eta}$ is the element of $\mathfrak b$ with 
$\exp^{\perp}_{\ast\xi}(\bar{\eta})=\eta$ (where $\bar{\eta}$ is regarded as 
an element of $T_{\xi}\mathfrak p'$ under the natural identification of 
$\mathfrak p'$ with $T_{\xi}\mathfrak p'$).  
Also, for each $Y\in T_x(M\cap\exp^{\perp}(\mathfrak p'))
\cap g_{\ast}\mathfrak p_{\beta}$ ($\beta\in\overline{\triangle}^V_+$), 
we can show 
$$A_{\eta}Y=-\frac{\beta(\bar{\eta})}{\tanh\beta(\xi)}Y
\leqno{({\rm A}.2)}$$
(see the proof of Theorem B of [Koi3]).  Let $\widetilde M_0$ be a component 
of $\widetilde M:=(\pi\circ\phi)^{-1}(M)$ and $\widetilde A$ be the shape 
tensor of $\widetilde M_0$.  From $({\rm A}.1),\,({\rm A}.2)$ and Lemma 5.1, 
we have 
$$\begin{array}{l}
\displaystyle{{\rm Spec}\widetilde A^{\bf c}_{\eta^L}=\{0\}\cup
\{\frac{\beta(\bar{\eta})}{-\beta(\xi)+j\pi\sqrt{-1}}\,\vert\,
\beta\in\overline{\triangle}_+^V,\,j\in{\bf Z}\}}\\
\hspace{2.3truecm}\displaystyle{\cup
\{\frac{\beta(\bar{\eta})}{-\beta(\xi)+(j+\frac12)\pi\sqrt{-1}}\,\vert\,
\beta\in\overline{\triangle}_+^H,\,j\in{\bf Z}\}.}
\end{array}
\leqno{({\rm A}.3)}$$
Denote by $F$ the sum of all complex focal hyperplanes of 
$(\widetilde M_0,u)$, where $u\in(\pi\circ\phi)^{-1}(gK)\cap\widetilde M_0$.  
From $({\rm A}.3)$, we have 
$$\begin{array}{l}
\displaystyle{g_{\ast}^{-1}F=\left(\mathop{\cup}_{\beta\in
\overline{\triangle}_+^V}
\mathop{\cup}_{j\in{\bf Z}}
(\beta^{\bf c})^{-1}(-\beta(\xi)+j\pi\sqrt{-1})\right)}\\
\hspace{1.7truecm}\displaystyle{\cup\left(
\mathop{\cup}_{\beta\in\overline{\triangle}_+^H}
\mathop{\cup}_{j\in{\bf Z}}(\beta^{\bf c})^{-1}
(-\beta(\xi)+(j+\frac12)\pi\sqrt{-1})\right)}\\
\hspace{1.2truecm}\displaystyle{=\left(\mathop{\cup}_{\beta\in
\overline{\triangle}_+^V}\mathop{\cup}_{j\in{\bf Z}}
(-\xi+(\beta^{\bf c})^{-1}(j\pi\sqrt{-1}))
\right)}\\
\hspace{1.7truecm}\displaystyle{\cup\left(
\mathop{\cup}_{\beta\in\overline{\triangle}_+^H}\mathop{\cup}_{j\in{\bf Z}}
(-\xi+(\beta^{\bf c})^{-1}((j+\frac12)\pi\sqrt{-1}))\right),}
\end{array}\leqno{({\rm A}.4)}$$
where we regard $F$ as a subspace of $(T^{\perp}_{gK}M)^{\bf c}$ under the 
natural identification of $(T^{\perp}_u\widetilde M_0)^{\bf c}$ with 
$(T^{\perp}_{gK}M)^{\bf c}$.  
From $({\rm A}.4)$, it follows that 
the complex Coxeter group $W_M$ associated 
with $M$ is isomorphic to the affine Weyl group associated with 
the root system $\overline{\triangle}$.  
Also, we have $g_{\ast}^{-1}(F\cap T^{\perp}_{gK}M)
=\displaystyle{\mathop{\cup}_{\beta\in\overline{\triangle}_+^V}
(-\xi+\beta^{-1}(0))}$.  Hence the real Coxeter group 
$W_{M,{\bf R}}$ associated with $M$ is isomorphic to the group generated by 
the reflections with respect to $\beta^{-1}(0)$'s ($\beta\in
\overline{\triangle}_+^V$).  Since 
$\overline{\triangle}^V_+$ is the positive root system associated with 
the symmetric space $\exp^{\perp}(\mathfrak p')=L/H\cap K$, 
$W_{M,{\bf R}}$ is isomorphic to the Weyl group associated with $L/H\cap K$.  

Next we shall list up the numbers 
$\displaystyle{\mathop{{\rm max}}_{v\in T^{\perp}M}\sharp{\rm Spec}\,A_v}$ 
for principal orbits $M$'s of Hermann type actions $H$'s on 
irreducible symmetric spaces $G/K$'s of non-compact type satisfying 
${\rm cohom}\,H={\rm rank}(G/K)$.  
We shall use the notations of the last paragraph.  
Since $\mathfrak p_{\beta}=\mathfrak p_{\beta}\cap\mathfrak p'+
\mathfrak p_{\beta}\cap{\mathfrak p'}^{\perp}$ for each 
$\beta\in\overline{\triangle}_+$, we have $\overline{\triangle}_+=
\overline{\triangle}_+^V\cup\overline{\triangle}_+^H$.  
Hence, from $(A.1)$ and $(A.2)$, we have the following equality:
$$\mathop{{\rm max}}_{v\in T^{\perp}M}\sharp{\rm Spec}\,A_v
=\sharp\overline{\triangle}_++\sharp(\overline{\triangle}^V_+\cap
\overline{\triangle}^H_+).\leqno{(A.5)}$$
In case of ${\rm cohom}\,H={\rm rank}(G/K)$, then we have 
$\mathfrak a=\mathfrak b$ and hence $\overline{\triangle}_+=\triangle_+$.  
Hence we can list up the numbers 
$\displaystyle{\mathop{{\rm max}}_{v\in T^{\perp}M}\sharp{\rm Spec}\,A_v}$ 
for the principal orbits $M$'s in the case (see Tables 2$\sim$4).  
The symbol $\widetilde{SO_0(1,8)}$ in Table 4 denotes 
the universal covering of $SO_0(1,8)$ and the symbol $\alpha$ in Table 4 
denotes an outer automorphism of $G_2^2$.  


\vspace{0.5truecm}

$$\begin{tabular}{|c|c|c|}
\hline
{\scriptsize$H$} & {\scriptsize$G/K$} & 
{\scriptsize$\mathop{{\rm max}}_{v\in T^{\perp}M}\sharp{\rm Spec}\,A_v$} \\
\hline
{\scriptsize$SO(n)$} & {\scriptsize$SL(n,{\bf R})/SO(n)$} & 
{\scriptsize$\frac{n(n-1)}{2}$}\\
\hline
{\scriptsize$SO_0(p,n-p)$} & {\scriptsize$SL(n,{\bf R})/SO(n)$} & 
{\scriptsize$\frac{n(n-1)}{2}$}\\
\hline
{\scriptsize$Sp(n)$} & {\scriptsize$SU^{\ast}(2n)/Sp(n)$} & 
{\scriptsize$\frac{n(n-1)}{2}$}\\
\hline
{\scriptsize$SO^{\ast}(2n)$} & {\scriptsize$SU^{\ast}(2n)/Sp(n)$} 
& {\scriptsize $n(n-1)$}\\
\hline
{\scriptsize$Sp(p,n-p)$} & {\scriptsize$SU^{\ast}(2n)/Sp(n)$} & 
{\scriptsize$\frac{n(n-1)}{2}$}\\
\hline
{\scriptsize$S(U(p)\times U(q))\,\,(p\leq q)$} & 
{\scriptsize$SU(p,q)/S(U(p)\times U(q))$} 
& {\scriptsize $\displaystyle{\left\{
\begin{array}{ll}
p(p+1) & (p<q)\\
p^2 & (p=q)
\end{array}\right.}$}\\
\hline
{\scriptsize$SO_0(p,q)\,\,(p\leq q)$} & 
{\scriptsize$SU(p,q)/S(U(p)\times U(q))$} & 
{\scriptsize$p(2p+1)$}\\
\hline
{\scriptsize$SO^{\ast}(2p)$} & {\scriptsize$SU(p,p)/S(U(p)\times U(p))$} & 
{\scriptsize$p(2p-1)$}\\
\hline
{\scriptsize$SL(p,{\bf C})\cdot U(1)$} 
& {\scriptsize$SU(p,p)/S(U(p)\times U(p))$} & {\scriptsize$p^2$}\\
\hline
{\scriptsize$SU(n)$} & {\scriptsize$SL(n,{\bf C})/SU(n)$} & 
{\scriptsize$\frac{n(n-1)}{2}$}\\
\hline
{\scriptsize$SO(n,{\bf C})$} & 
{\scriptsize$SL(n,{\bf C})/SU(n)$} & {\scriptsize$n(n-1)$}\\
\hline
{\scriptsize$SO(p)\times SO(q)\,\,(p\leq q)$} & 
{\scriptsize$SO_0(p,q)/SO(p)\times SO(q)$} 
& {\scriptsize $\displaystyle{\left\{
\begin{array}{ll}
p^2 & (p<q)\\
p(p-1) & (p=q)
\end{array}\right.}$}\\
\hline
{\scriptsize$SO(p,{\bf C})$} & {\scriptsize$SO_0(p,p)/SO(p)\times SO(p)$} & 
{\scriptsize$p(p-1)$}\\
\hline
{\scriptsize$U(n)$} & {\scriptsize$SO^{\ast}(2n)/U(n)$} & 
{\scriptsize$\displaystyle{\left\{
\begin{array}{ll}
\frac{n^2-1}{4} & (n:{\rm odd})\\
\frac{n^2}{4} & (n:{\rm even})
\end{array}\right.}$}\\
\hline
{\scriptsize$SO(n,{\bf C})$} & {\scriptsize$SO^{\ast}(2n)/U(n)$} 
& {\scriptsize$\frac{n(n-1)}{2}$}\\
\hline
{\scriptsize$SU(2i,2n-2i)\cdot U(1)$} & {\scriptsize$SO^{\ast}(4n)/U(2n)$} 
& {\scriptsize$n^2$}\\
\hline
{\scriptsize$SU(i,2n-i+1)\cdot U(1)$} & {\scriptsize$SO^{\ast}(4n+2)/U(2n+1)$} 
& {\scriptsize$n^2+n$}\\
\hline
{\scriptsize$SO_0(i,2n-i+1)$} & 
{\scriptsize$SO(2n+1,{\bf C})/SO(2n+1)$} & {\scriptsize$2n^2$}\\
\hline
{\scriptsize$SO_0(2i,2n-2i)$} & 
{\scriptsize$SO(2n,{\bf C})/SO(2n)$} & {\scriptsize$\frac{(2n-1)^2}{2}$}\\
\hline
{\scriptsize$U(n)$} & {\scriptsize$Sp(n,{\bf R})/U(n)$} & {\scriptsize$n^2$}\\
\hline
{\scriptsize$SU(i,n-i)\cdot U(1)$} & {\scriptsize$Sp(n,{\bf R})/U(n)$} & 
{\scriptsize$n^2$}\\
\hline
{\scriptsize$Sp(p)\times Sp(q)$} & 
{\scriptsize$Sp(p,q)/Sp(p)\times Sp(q)$} & 
{\scriptsize$\displaystyle{\left\{
\begin{array}{ll}
p^2+p & (p<q)\\
p^2 & (p=q)
\end{array}\right.}$}\\
\hline
{\scriptsize$SU(p,q)\cdot U(1)$} & {\scriptsize$Sp(p,q)/Sp(p)\times Sp(q)$} & 
{\scriptsize $\frac12p(3p+5)$}\\
\hline
{\scriptsize$SU^{\ast}(2p)\cdot U(1)$} & 
{\scriptsize$Sp(p,p)/Sp(p)\times Sp(p)$} & {\scriptsize$2p^2$}\\
\hline
{\scriptsize$SL(n,{\bf C})\cdot SO(2,{\bf C})$} 
& {\scriptsize$Sp(n,{\bf C})/Sp(n)$} & {\scriptsize$2n^2$}\\
\hline
{\scriptsize$Sp(n,{\bf R})$} & {\scriptsize$Sp(n,{\bf C})/Sp(n)$} 
& {\scriptsize$n^2$}\\
\hline
{\scriptsize$Sp(i,n-i)$} & {\scriptsize$Sp(n,{\bf C})/Sp(n)$} & 
{\scriptsize$n^2$}\\
\hline
\end{tabular}$$

\centerline{{\bf Table 2.}}

\newpage

$$\begin{tabular}{|c|c|c|}
\hline
{\scriptsize$H$} & {\scriptsize$G/K$} & 
{\scriptsize$\mathop{{\rm max}}_{v\in T^{\perp}M}\sharp{\rm Spec}\,A_v$} \\
\hline
{\scriptsize$Sp(4)/\{\pm1\}$} & {\scriptsize$E_6^6/(Sp(4)/\{\pm1\})$} 
& {\scriptsize$36$}\\
\hline
{\scriptsize$Sp(4,{\bf R})$} & {\scriptsize$E_6^6/(Sp(4)/\{\pm1\})$} & 
{\scriptsize$36$}\\
\hline
{\scriptsize$Sp(2,2)$} & {\scriptsize$E_6^6/(Sp(4)/\{\pm1\})$} & 
{\scriptsize$36$}\\
\hline
{\scriptsize$SU(6)\cdot SU(2)$} & {\scriptsize$E_6^2/SU(6)\cdot SU(2)$} & 
{\scriptsize$24$}\\
\hline
{\scriptsize$Sp(1,3)$} & {\scriptsize$E_6^2/SU(6)\cdot SU(2)$} & 
{\scriptsize$36$}\\
\hline
{\scriptsize$Sp(4,{\bf R})$} & {\scriptsize$E_6^2/SU(6)\cdot SU(2)$} & 
{\scriptsize$34$}\\
\hline
{\scriptsize$SU(2,4)\cdot SU(2)$} & {\scriptsize$E_6^2/SU(6)\cdot SU(2)$} & 
{\scriptsize$30$}\\
\hline
{\scriptsize$SU(3,3)\cdot SL(2,{\bf R})$} 
& {\scriptsize$E_6^2/SU(6)\cdot SU(2)$} & {\scriptsize$24$}\\
\hline
{\scriptsize$Spin(10)\cdot U(1)$} & 
{\scriptsize$E_6^{-14}/Spin(10)\cdot U(1)$} & {\scriptsize$6$}\\
\hline
{\scriptsize$Sp(2,2)$} & 
{\scriptsize$E_6^{-14}/Spin(10)\cdot U(1)$} & {\scriptsize$10$}\\
\hline
{\scriptsize$SU(2,4)\cdot SU(2)$} 
& {\scriptsize$E_6^{-14}/Spin(10)\cdot U(1)$} & {\scriptsize$10$}\\
\hline
{\scriptsize$SU(1,5)\cdot SL(2,{\bf R})$} 
& {\scriptsize$E_6^{-14}/Spin(10)\cdot U(1)$} & {\scriptsize$10$}\\
\hline
{\scriptsize$SO^{\ast}(10)\cdot U(1)$} 
& {\scriptsize$E_6^{-14}/Spin(10)\cdot U(1)$} & {\scriptsize$7$}\\
\hline
{\scriptsize$SO_0(2,8)\cdot U(1)$} & 
{\scriptsize$E_6^{-14}/Spin(10)\cdot U(1)$} & {\scriptsize$10$}\\
\hline
{\scriptsize$F_4$} & {\scriptsize$E_6^{-26}/F_4$} & {\scriptsize$3$}\\
\hline
{\scriptsize$Sp(1,3)$} & {\scriptsize$E_6^{-26}/F_4$} & {\scriptsize$6$}\\
\hline
{\scriptsize$F_4^{-20}$} & {\scriptsize$E_6^{-26}/F_4$} & {\scriptsize$3$}\\
\hline
{\scriptsize$E_6$} & {\scriptsize$E_6^{\bf c}/E_6$} & {\scriptsize$36$}\\
\hline
{\scriptsize$E_6^2$} & {\scriptsize$E_6^{\bf c}/E_6$} & {\scriptsize$36$}\\
\hline
{\scriptsize$E_6^{-14}$} & {\scriptsize$E_6^{\bf c}/E_6$} & {\scriptsize$36$}\\
\hline
{\scriptsize$Sp(4,{\bf C})$} & {\scriptsize$E_6^{\bf c}/E_6$} & 
{\scriptsize$72$}\\
\hline
{\scriptsize$SU(8)/\{\pm1\}$} & {\scriptsize$E_7^7/(SU(8)/\{\pm1\})$} & 
{\scriptsize$63$}\\
\hline
{\scriptsize$SL(8,{\bf R})$} & {\scriptsize$E_7^7/(SU(8)/\{\pm1\})$} & 
{\scriptsize$63$}\\
\hline
{\scriptsize$SU^{\ast}(8)$} & {\scriptsize$E_7^7/(SU(8)/\{\pm1\})$} & 
{\scriptsize$63$}\\
\hline
{\scriptsize$SU(4,4)$} & {\scriptsize$E_7^7/(SU(8)/\{\pm1\})$} & 
{\scriptsize$63$}\\
\hline
{\scriptsize$SO'(12)\cdot SU(2)$} & 
{\scriptsize$E_7^{-5}/SO'(12)\cdot SU(2)$} & 
{\scriptsize$24$}\\
\hline
{\scriptsize$SU(4,4)$} & {\scriptsize$E_7^{-5}/SO'(12)\cdot SU(2)$} & 
{\scriptsize$36$}\\
\hline
{\scriptsize$SU(2,6)$} & {\scriptsize$E_7^{-5}/SO'(12)\cdot SU(2)$} & 
{\scriptsize$36$}\\
\hline
{\scriptsize$SO^{\ast}(12)\cdot SL(2,{\bf R})$} & 
{\scriptsize$E_7^{-5}/SO'(12)\cdot SU(2)$} & 
{\scriptsize$24$}\\
\hline
{\scriptsize$SO_0(4,8)\cdot SU(2)$} & 
{\scriptsize$E_7^{-5}/SO'(12)\cdot SU(2)$} & {\scriptsize$24$}\\
\hline
{\scriptsize$E_6\cdot U(1)$} & 
{\scriptsize$E_7^{-25}/E_6\cdot U(1)$} & 
{\scriptsize$9$}\\
\hline
{\scriptsize$SU^{\ast}(8)$} & {\scriptsize$E_7^{-25}/E_6\cdot U(1)$} & 
{\scriptsize$15$}\\
\hline
{\scriptsize$SU(2,6)$} & {\scriptsize$E_7^{-25}/E_6\cdot U(1)$} & 
{\scriptsize$15$}\\
\hline
{\scriptsize$E_6^{-14}\cdot U(1)$} & 
{\scriptsize$E_7^{-25}/E_6\cdot U(1)$} & {\scriptsize$9$}\\
\hline
{\scriptsize$E_7$} & {\scriptsize$E_7^{\bf c}/E_7$} & {\scriptsize$63$}\\
\hline
{\scriptsize$E_7^7$} & {\scriptsize$E_7^{\bf c}/E_7$} & {\scriptsize$63$}\\
\hline
{\scriptsize$E_7^{-5}$} & {\scriptsize$E_7^{\bf c}/E_7$} & {\scriptsize$63$}\\
\hline
{\scriptsize$E_7^{-25}$} & {\scriptsize$E_7^{\bf c}/E_7$} & {\scriptsize$63$}\\
\hline
{\scriptsize$SL(8,{\bf C})$} & {\scriptsize$E_7^{\bf c}/E_7$} & 
{\scriptsize$126$}\\
\hline
\end{tabular}$$

\centerline{{\bf Table 3.}}


\vspace{1truecm}

$$\begin{tabular}{|c|c|c|}
\hline
{\scriptsize$H$} & {\scriptsize$G/K$} & 
{\scriptsize$\mathop{{\rm max}}_{v\in T^{\perp}M}\sharp{\rm Spec}\,A_v$} \\
\hline
{\scriptsize$SO'(16)$} & {\scriptsize$E_8^8/SO'(16)$} & {\scriptsize$120$}\\
\hline
{\scriptsize$SO_0(8,8)$} & {\scriptsize$E_8^8/SO'(16)$} & 
{\scriptsize$120$}\\
\hline
{\scriptsize$E_7\cdot Sp(1)$} & {\scriptsize$E_8^{-24}/E_7\cdot Sp(1)$} & 
{\scriptsize$24$}\\
\hline
{\scriptsize$SO^{\ast}(16)$} & {\scriptsize$E_8^{-24}/E_7\cdot Sp(1)$} & 
{\scriptsize$36$}\\
\hline
{\scriptsize$SO_0(4,12)$} & {\scriptsize$E_8^{-24}/E_7\cdot Sp(1)$} & 
{\scriptsize$36$}\\
\hline
{\scriptsize$E_7^{-5}\cdot Sp(1)$} & {\scriptsize$E_8^{-24}/E_7\cdot Sp(1)$} & 
{\scriptsize$24$}\\
\hline
{\scriptsize$E_7^{-25}\cdot SL(2,{\bf R})$} & 
{\scriptsize$E_8^{-24}/E_7\cdot Sp(1)$} & {\scriptsize$24$}\\
\hline
{\scriptsize$E_8$} & {\scriptsize$E_8^{\bf c}/E_8$} & {\scriptsize$120$}\\
\hline
{\scriptsize$E_8^8$} & {\scriptsize$E_8^{\bf c}/E_8$} & {\scriptsize$120$}\\
\hline
{\scriptsize$E_8^{-24}$} & {\scriptsize$E_8^{\bf c}/E_8$} & 
{\scriptsize$120$}\\
\hline
{\scriptsize$SO(16,{\bf C})$} & {\scriptsize$E_8^{\bf c}/E_8$} & 
{\scriptsize$240$}\\
\hline
{\scriptsize$Sp(3)\cdot Sp(1)$} & {\scriptsize$F_4^4/Sp(3)\cdot Sp(1)$} & 
{\scriptsize$24$}\\
\hline
{\scriptsize$Sp(1,2)\cdot Sp(1)$} & {\scriptsize$F_4^4/Sp(3)\cdot Sp(1)$} 
& {\scriptsize$24$}\\
\hline
{\scriptsize$Sp(3,{\bf R})\cdot SL(2,{\bf R})$} 
& {\scriptsize$F_4^4/Sp(3)\cdot Sp(1)$} & {\scriptsize$24$}\\
\hline
{\scriptsize$Spin(9)$} & {\scriptsize$F_4^{-20}/Spin(9)$} 
& {\scriptsize$2$}\\
\hline
{\scriptsize$Sp(1,2)\cdot Sp(1)$} & {\scriptsize$F^{-20}_4/Spin(9)$} & 
{\scriptsize$2$}\\
\hline
{\scriptsize$\widetilde{SO_0(1,8)}$} & {\scriptsize$F^{-20}_4/Spin(9)$} & 
{\scriptsize$4$}\\
\hline
{\scriptsize$F_4$} & {\scriptsize$F_4^{\bf C}/F_4$} & 
{\scriptsize$24$} \\
\hline
{\scriptsize$F_4^4$} & {\scriptsize$F_4^{\bf C}/F_4$} & 
{\scriptsize$24$} \\
\hline
{\scriptsize$F_4^{-20}$} & {\scriptsize$F_4^{\bf C}/F_4$} & {\scriptsize$24$}\\
\hline
{\scriptsize$Sp(3,{\bf C})\cdot SL(2,{\bf C})$} 
& {\scriptsize$F_4^{\bf C}/F_4$} & {\scriptsize$48$}\\
\hline
{\scriptsize$SO(4)$} & {\scriptsize$G_2^2/SO(4)$} & {\scriptsize$6$}\\
\hline
{\scriptsize$SL(2,{\bf R})\times SL(2,{\bf R})$} & {\scriptsize$G_2^2/SO(4)$} 
& {\scriptsize$6$}\\
\hline
{\scriptsize$\alpha(SO(4))$} & {\scriptsize$G_2^2/SO(4)$} & 
{\scriptsize$6$}\\
\hline
{\scriptsize$G_2$} & {\scriptsize$G_2^{\bf c}/G_2$} & {\scriptsize$6$}\\
\hline
{\scriptsize$G_2^2$} & {\scriptsize$G_2^{\bf c}/G_2$} & {\scriptsize$6$}\\
\hline
{\scriptsize$SL(2,{\bf C})\times SL(2,{\bf C})$} 
& {\scriptsize$G_2^{\bf c}/G_2$} & {\scriptsize$12$}\\
\hline
\end{tabular}$$

\centerline{{\bf Table 4.}}

\vspace{1truecm}

\centerline{{\bf Appendix 2}}

\vspace{0.5truecm}

In this appendix, we prove the following important fact for 
a curvature-adapted submanifold with globally flat and abelian normal bundle 
in a symmetric space.  

\vspace{0.5truecm}

\noindent
{\bf Proposition A.1.} {\sl Let $M$ be a curvature-adapted submanifold 
with globally flat and abelian normal bundle in a symmetric space $G/K$, $A$ 
be the shape tensor of $M$ and $R$ be the curvature tensor of $G/K$.  
Then, for any $x\in M$, 
$$\{R(\cdot,v)v\vert_{T_xM}\,\vert\,v\in T_x^{\perp}M\}\cup
\{A_v\,\vert\,v\in T^{\perp}_xM\}$$
is a commuting family of linear transformations of $T_xM$.}

\vspace{0.5truecm}

\noindent
{\it Proof.} We shall show this statement in the case where $G/K$ is of 
non-compact type.  Let $v_i\in T^{\perp}_xM$ ($i=1,2$).  Since $M$ has abelian 
normal bundle, $R(\cdot,v_1)v_1\vert_{T_xM}$ and $R(\cdot,v_2)v_2\vert_{T_xM}$ 
commute with each other.  Since $M$ has globally flat and abelian normal 
bundle, $A_{v_1}$ and $A_{v_2}$ commute with each other.  We shall show that 
$R(\cdot,v_1)v_1\vert_{T_xM}$ and $A_{v_2}$ commute with each other.  Let 
$x=gK$.  Take a maximal abelian subspace $\mathfrak a$ of 
$\mathfrak p:=T_{eK}(G/K)$ containing $\mathfrak b:=g_{\ast}^{-1}
(T^{\perp}_xM)$.  Let $\triangle$ be the root system with respect to 
$\mathfrak a$ and set $\overline{\triangle}:=\{\alpha\vert_{\mathfrak b}
\,\vert\,\alpha\in\triangle\,\,{\rm s.t.}\,\,
\alpha\vert_{\mathfrak b}\not=0\}$.  For each $\beta\in\overline{\triangle}$, 
we set $\mathfrak p_{\beta}:=\{X\in\mathfrak p\,\vert\,{\rm ad}(b)^2(X)
=\beta(b)^2X\,\,(\forall\,b\in\mathfrak b)\}$.  Then we have 
$\mathfrak p=\mathfrak z_{\mathfrak p}(\mathfrak b)
+\sum_{\beta\in\overline{\triangle}_+}\mathfrak p_{\beta}$, where 
$\overline{\triangle}_+$ is the positive root system under some lexicographic 
ordering and $\mathfrak z_{\mathfrak p}(\mathfrak b)$ is the centralizer of 
$\mathfrak b$ in $\mathfrak p$.  Consider 
$$D:=\{v\in T^{\perp}_xM\,\vert\,{\rm Span}\{v\}\cap\left(
\mathop{\cup}_{(\beta_1,\beta_2)\in\overline{\triangle}_+\times\overline{\triangle}_+\,\,{\rm s.t.}\,\,\beta_1\not=\beta_2}({\it l}_{\beta_1}\cap
{\it l}_{\beta_2})\right)=\emptyset\},$$
where ${\it l}_{\beta_i}:=\beta_i^{-1}(1)$ ($i=1,2$).  
It is clear that $D$ is open and dense in $T^{\perp}_xM$.  
Take $v\in D$.  Then, since $\beta(v)$'s ($\beta\in\overline{\triangle}_+$) 
are mutually distinct, the decomposition $T_xM=g_{\ast}
(\mathfrak z_{\mathfrak p}(\mathfrak b)\ominus\mathfrak b)
+\sum_{\beta\in\overline{\triangle}_+}g_{\ast}\mathfrak p_{\beta}$ is the 
eigenspace decomposition of $R(\cdot,v)v\vert_{T_xM}$.  Since $M$ is 
curvature-adapted by the assumption and hence 
$[R(\cdot,v)v\vert_{T_xM},A_v]=0$, we have 
$$T_xM=\sum_{\lambda\in{\rm Spec}\,A_v}\left(
(g_{\ast}(\mathfrak z_{\mathfrak p}(\mathfrak b)\ominus\mathfrak b)\cap
{\rm Ker}(A_v-\lambda{\rm id}))+\sum_{\beta\in\overline{\triangle}_+}
(g_{\ast}\mathfrak p_{\beta}\cap{\rm Ker}(A_v-\lambda\,{\rm id}))\right).
\leqno{(A.6)}$$
Suppose that (A.6) does not hold for some $v_0\in T^{\perp}_xM\setminus D$.  
Then it is easy to show that there exists a neighborhood $U$ of $v_0$ in 
$T^{\perp}_xM$ such that (A.6) does not hold for any $v\in U$.  Clearly we 
have $U\cap D=\emptyset$.  This contradicts the fact that $D$ is dense in 
$T^{\perp}_xM$.  Hence (A.6) holds for any $v\in T^{\perp}_xM\setminus D$.  
Therefore, (A.6) holds for any $v\in T^{\perp}_xM$.  In particular, (A.6) 
holds for $v_2$.  On the other hand, the decomposition $T_xM=g_{\ast}
\mathfrak z_{\mathfrak p}(\mathfrak b)+\sum_{\beta\in\overline{\triangle}_+}
g_{\ast}\mathfrak p_{\beta}$ is the common eigenspace decomposition of 
$R(\cdot,v)v\vert_{T_xM}$'s ($v\in T^{\perp}_xM$).  From these facts, we have 
$$T_xM=\sum_{\lambda\in{\rm Spec}\,A_{v_2}}
\sum_{\mu\in{\rm Spec}\,R(\cdot,v_1)v_1\vert_{T_xM}}
\left({\rm Ker}(R(\cdot,v_1)v_1\vert_{T_xM}-\mu\,{\rm id})
\cap{\rm Ker}(A_{v_2}-\lambda\,{\rm id})\right),$$
which implies that $R(\cdot,v_1)v_1\vert_{T_xM}$ and $A_{v_2}$ commute with 
each other.  This completes the proof.  \hspace{12.9truecm}q.e.d.

\vspace{0.5truecm}

\noindent
{\it Remark A.1.} O. Goertsches and G. Thorbergsson [GT] have already shown 
that the statement of this proposition holds for principal orbits of Heremann 
actions on symmetric spaces of compact type.  



\vspace{1truecm}

\centerline{{\bf References}}

\vspace{0.5truecm}

{\small

\noindent
[B1] J. Berndt, Real hypersurfaces with constant principal curvatures 
in complex hyperbolic 

space, J. Reine Angew. Math. {\bf 395} (1989) 132-141. 

\noindent
[B2] J. Berndt, Real hypersurfaces in quaternionic space forms, 
J. Reine Angew. Math. {\bf 419} 

(1991) 9-26. 





\noindent
[BT] J. Berndt and H. Tamaru, Cohomogeneity one actions on noncompact 
symmetric spaces 

with a totally geodesic singular orbit, 
Tohoku Math. J. {\bf 56} (2004) 163-177.

\noindent
{\small [BV] J. Berndt and L. Vanhecke, 
Curvature adapted submanifolds, 
Nihonkai Math. J. {\bf 3} (1992) 

177-185.



\noindent
[Ch] U. Christ, 
Homogeneity of equifocal submanifolds, J. Differential Geometry 
{\bf 62} (2002) 1-15.

\noindent
[Co] H.S.M. Coxeter, 
Discrete groups generated by reflections, 
Ann. of Math. {\bf 35} (1934) 588-621.

\noindent
[E] H. Ewert, A splitting theorem for equifocal submanifolds in simply 
connected compact symme-

tric spaces, Proc. of Amer. Math. Soc. {\bf 126} (1998) 2443-2452.


\noindent
[G1] L. Geatti, 
Invariant domains in the complexfication of a noncompact Riemannian 
symmetric 

space, J. of Algebra {\bf 251} (2002) 619-685.

\noindent
[G2] L. Geatti, 
Complex extensions of semisimple symmetric spaces, manuscripta math. {\bf 120} 

(2006) 1-25.

\noindent
[GT] O. Goertsches and G. Thorbergsson, 
On the Geometry of the orbits of Hermann actions, 

arXiv:math.DG/0701731.



\noindent
[HLO] E. Heintze, X. Liu and C. Olmos, 
Isoparametric submanifolds and a 
Chevalley type rest-

riction theorem, Integrable systems, geometry, and topology, 151-190, 
AMS/IP Stud. Adv. 

Math. 36, Amer. Math. Soc., Providence, RI, 2006.

\noindent
[HPTT] E. Heintze, R.S. Palais, C.L. Terng and G. Thorbergsson, 
Hyperpolar actions on symme-

tric spaces, Geometry, topology and physics for Raoul Bott (ed. S. T. Yau), Conf. Proc. 

Lecture Notes Geom. Topology {\bf 4}, 
Internat. Press, Cambridge, MA, 1995 pp214-245.

\noindent
[He] S. Helgason, 
Differential geometry, Lie groups and symmetric spaces, 
Academic Press, New 

York, 1978.

\noindent
[Hu] M.C. Hughes, 
Complex reflection groups, 
Communications in Algebra {\bf 18} 
(1990) 3999-4029.

\noindent
[Ka] R. Kane, Reflection groups and Invariant Theory, CMS Books in 
Mathematics, Springer-

Verlag, New York, 2001.

\noindent
[Koi1] N. Koike, 
Submanifold geometries in a symmetric space of non-compact 
type and a pseudo-

Hilbert space, Kyushu J. Math. {\bf 58} (2004) 167-202.

\noindent
[Koi2] N. Koike, 
Complex equifocal submanifolds and infinite dimensional anti-
Kaehlerian isopara-

metric submanifolds, Tokyo J. Math. {\bf 28} (2005) 201-247.

\noindent
[Koi3] N. Koike, 
Actions of Hermann type and proper complex equifocal submanifolds, 
Osaka J. 

Math. {\bf 42} (2005) 599-611.

\noindent
[Koi4] N. Koike, 
A splitting theorem for proper complex equifocal submanifolds, Tohoku Math. J. 

{\bf 58} (2006) 393-417.



\noindent
[Koi5] N. Koike, A Chevalley type restriction theorem for a proper complex 
equifocal submani-

fold, Kodai Math. J. {\bf 30} (2007) 280-296.

\noindent
[Koi6] N. Koike, The complexifications of pseudo-Riemannian manifolds and 
anti-Kaehler geo-

metry, in preparation.

\noindent
[Kol] A. Kollross, A Classification of hyperpolar and cohomogeneity one 
actions, Trans. Amer. 

Math. Soc. {\bf 354} (2001) 571-612.




\noindent
[OS] T. Oshima and J. Sekiguchi, The restricted root system of a semisimple 
symmetric pair, 

Advanced Studies in Pure Math. {\bf 4} (1984), 433--497. 


\noindent
[PT] R.S. Palais and C.L. Terng, Critical point theory and submanifold 
geometry, Lecture Notes 

in Math. {\bf 1353}, Springer, Berlin, 1988.

\noindent
[R] W. Rossmann, 
The structure of semisimple symmetric spaces, Can. J. Math. {\bf 1} 
(1979), 157--180.

\noindent
[S1] R. Sz$\ddot{{{\rm o}}}$ke, Complex structures on tangent 
bundles of Riemannian manifolds, Math. Ann. {\bf 291} 

(1991) 409-428.

\noindent
[S2] R. Sz$\ddot{{{\rm o}}}$ke, Automorphisms of certain Stein 
manifolds, Math. Z. {\bf 219} (1995) 357-385.

\noindent
[S3] R. Sz$\ddot{{{\rm o}}}$ke, Adapted complex structures and 
geometric quantization, Nagoya Math. J. {\bf 154} 

(1999) 171-183.

\noindent
[S4] R. Sz$\ddot{{{\rm o}}}$ke, Involutive structures on the 
tangent bundle of symmetric spaces, 
Math. Ann. {\bf 319} 

(2001), 319--348.

\noindent
[S5] R. Sz$\ddot{{{\rm o}}}$ke, Canonical complex structures associated to 
connections and complexifications of 

Lie groups, Math. Ann. {\bf 329} (2004), 553--591.

\noindent
[T] C.L. Terng, 
Isoparametric submanifolds and their Coxeter groups, 
J. Differential Geometry 

{\bf 21} (1985) 79-107.




\noindent
[TT] C.L. Terng and G. Thorbergsson, 
Submanifold geometry in symmetric spaces, 
J. Differential 

Geometry {\bf 42} (1995) 665-718.


}

\vspace{1truecm}

\rightline{Department of Mathematics, Faculty of Science, }
\rightline{Tokyo University of Science}
\rightline{26 Wakamiya Shinjuku-ku,}
\rightline{Tokyo 162-8601, Japan}
\rightline{(e-mail: koike@ma.kagu.tus.ac.jp)}

\end{document}